\documentclass{irmaart}

\usepackage{amsbsy,amscd,amsfonts,latexsym,amstext,delarray, diagrams}

\numberwithin{equation}{section}

\theoremstyle{plain}
\newtheorem{theorem}{Theorem}[section]

\newtheorem{proposition}[theorem]{Proposition}

\newtheorem{definition}[theorem]{Definition}

\newcommand{\scr}{\mathcal}
\newcommand{\Hom}{\operatorname{Hom}}
\newcommand{\Ext}{{\rm Ext}}
\newcommand{\Ker}{\operatorname{Ker}}
\newcommand{\Spec}{{\rm Spec}}
\newcommand{\Image}{{\rm Im}}
\newcommand{\Corr}{{\rm Corr}}
\newcommand{\Gal}{{\rm Gal}}

\newcommand{\cf}{{\it cf.\/}\ }
\newcommand{\ie}{{\it i.e.\/}\ }
\newcommand{\eg}{{\it e.g.\/}\ }
\newcommand{\resp}{{\it resp.\/}\ }
\newcommand{\op}{{\it op.cit.\/}\ }
\newcommand{\no}{\noindent}

\def\qqq{\,,\quad \forall\,}

\newcommand{\N}{\mathbb{N}}
\newcommand{\Q}{\mathbb{Q}}
\newcommand{\Z}{\mathbb{Z}}
\newcommand{\C}{\mathbb{C}}
\newcommand{\R}{\mathbb{R}}

\newcommand{\bF}{{\bf F}}
\newcommand{\bG}{\mathbb{G}}
\newcommand{\bP}{{\bf P}}
\newcommand{\cA}{\scr{A}}
\newcommand{\cB}{\scr{B}}
\newcommand{\cC}{\scr{C}}
\newcommand{\cD}{\scr{D}}
\newcommand{\cE}{\scr{E}}
\newcommand{\cF}{\scr{F}}
\newcommand{\cG}{\scr{G}}
\newcommand{\cK}{\scr{K}}
\newcommand{\cM}{\scr{M}}
\newcommand{\cO}{\scr{O}}
\newcommand{\cV}{\scr{V}}
\newcommand{\cZ}{\scr{Z}}

\begin{document}

\address{
Mathematics Department, The Johns Hopkins
University \\ 3400 N. Charles Street, Baltimore MD 21218 USA\\
email:\,\tt{kc@math.jhu.edu}
}

\title{Noncommutative geometry and motives
\\(a quoi servent les endomotifs?)}

\author{Caterina Consani\thanks{Work partially supported by NSF-FRG grant DMS-0652431. The author
wishes to thank the organizers of the
Conference for their kind invitation to speak and the CIRM in
Luminy-Marseille for the pleasant atmosphere and their support.}}

\maketitle

\begin{abstract}
This paper gives a short and historical survey on the theory of pure motives
in algebraic geometry and reviews some of the recent developments of this
theory in noncommutative geometry. The second part of the paper outlines
the new theory of endomotives and some of its relevant 
applications in number-theory. 
\end{abstract}

\section{Introduction}

This paper is based on  three lectures I gave at the Conference on
``Renormalization and Galois theories'' that was held in Luminy, at
the Centre International de Rencontres Math\'ematiques (CIRM), on
March 2006. The purpose of these talks was to give an elementary
overview on classical motives (pure motives) and to survey on
some of the recent developments of this theory in noncommutative
geometry, especially following the introduction of the notion of an
{\it endomotive}.

It is likely to expect that the reader acquainted with the literature
on motives theory will not fail to notice the allusion, in this title,
to the  paper \cite{17} in which P. Deligne states that in spite of the lack of
essential progresses on the problem of constructing 
``relevant'' algebraic cycles, the techniques supplied by the
theory of motives remain a powerful tool in algebraic geometry and
arithmetic.

\no The assertion on the lack of relevant progresses on algebraic
cycles seems, unfortunately, still to apply at the present time, fifteen years
after Deligne wrote his paper. Despite the general failure of testing 
the Standard Conjectures, it is also
true that in these recent years the knowledge on motives has been
substantially improved by several new results and also by some unexpected
developments.

Motives were introduced by A. Grothendieck with the aim to supply
an intrinsic explanation for the analogies occurring among various
cohomological theories in algebraic geometry. They are
expected to play the role of a universal cohomological theory by also
furnishing a linearization of the theory of  algebraic varieties and in the
original understanding they were expected to provide the correct framework for a successful 
approach to the Weil's Conjectures on the zeta-function of a variety over a finite field.

\no Even though the Weil's Conjectures have been proved by Deligne
without appealing to the theory of motives, an enlarged and in part still conjectural
theory of mixed motives has in the meanwhile proved its usefulness
in explaining conceptually, some intriguing phenomena arising in
several areas of pure mathematics, such as Hodge theory, $K$-theory,
algebraic cycles, polylogarithms, $L$-functions, Galois
representations etc.

Very recently, some new developments of the theory of motives to
number-theory and quantum field theory  have been found or are about
to be developed, with the support of techniques supplied by
noncommutative geometry and the theory of operator algebras.

\no In number-theory, a conceptual understanding of the main result of
\cite{10} on the interpretation proposed by A. Connes of the Weil
explicit formulae as a Lefschetz trace formula over the
noncommutative space of ad\`ele classes, requires the introduction
of a generalized category of motives  inclusive of spaces
which are highly singular from a classical viewpoint.

\no The problem of finding a suitable enlargement of the category of
(smooth projective) algebraic varieties is combined with the even
more compelling one of the definition of a generalized notion of
correspondences. Several  questions arise already when one considers
special types of zero-dimensional noncommutative spaces, such as the
space underlying the quantum statistical dynamical system defined by
J. B. Bost and Connes in \cite{6} (the BC-system). This space
is a simplified version of the ad\`eles class space of \cite{10} and
it encodes in its group of symmetries, the arithmetic of the maximal
abelian extension of $\Q$.

In this paper I give an overview on the theory of endomotives
(algebraic and analytic). This theory has been originally developed
in the joint paper \cite{7} with A. Connes and M. Marcolli and
has been applied already in our subsequent work \cite{11}. The category of
endomotives is the minimal one that makes it possible to
understand conceptually the role played by the absolute Galois group in
several dynamical systems that have been recently introduced in noncommutative
geometry as generalizations of the BC-system, which was our motivating and prototype example.

\no The category of endomotives is  a natural enlargement of the
category of Artin motives: the objects are noncommutative spaces
defined by semigroup actions on projective limits of Artin
motives. The morphisms generalize the notion of algebraic
correspondences  and are defined by means of \'etale groupoids to
account for the presence of the semigroup actions.

\no Endomotives carry a natural Galois action which is inherited from
the Artin motives and they have both an algebraic and an analytic
description. The latter is particularly useful as it provides the
data of a quantum statistical dynamical system, via the implementation of
a {\it canonical time evolution} (a one-parameter family of automorphisms) which is associated
by the theory of M.~Tomita (\cf\cite{36}) to an initial state (probability
measure) assigned on an analytic endomotive.
This is the crucial new development supplied by the theory of operator-algebras
to a noncommutative $C^*$-algebra and in particular to the algebra of the
BC-system. 

\no The implication in number-theory is striking: the time
evolution implements on the dual system a scaling action which 
combines with the action of the Galois group to determine on the
cyclic homology of a suitable noncommutative motive associated to
the original endomotive, a characteristic zero analog of the action
of the Weil group on the \'etale cohomology of an algebraic variety.
When these techniques are applied to the endomotive of the BC-system
or to the endomotive of the ad\`eles class space, the main 
implication is the spectral realization of the zeroes of the
corresponding $L$-functions.

These results supply a first answer to the question I raised in the
title of this paper (a quoi servent les endomotifs?). An open and interesting
problem is connected to the definition of a higher dimensional
theory of noncommutative motives and in particular the introduction of a
theory of noncommutative elliptic motives and modular forms. A related
problem is of course connected to the definition of a higher
dimensional theory of geometric correspondences. The comparison
between algebraic correspondences for motives and geometric
correspondences for noncommutative spaces is particularly easy in
the zero-dimensional case, because the equivalence relations play no
role. In noncommutative geometry, algebraic cycles are naturally
replaced by bi-modules, or by classes in equivariant $KK$-theory. Naturally, 
the original problem of finding ``interesting''
cycles pops-up again in this  topological framework:  a
satisfactory solution to this question seems to be one of the main steps to
undertake for a further development of these ideas.

\section{Classical motives: an overview}

The theory of motives in algebraic geometry was  established by A.
Grothendieck in the 1960s: 1963-69 (\cf\cite{35},\cite{21}). The
foundations are documented in the unpublished manuscript \cite{19}
and were discussed in a seminar at the Institut des Hautes \'Etudes
Scientifiques, in 1967. This theory was conceived 
as a fundamental machine to develop Grothendieck's ``long-run
program'' focused on the theme of the connections between geometry
and arithmetic. 

\no At the heart of the philosophy of motives sit
Grothendieck's speculations on the existence of a universal
cohomological theory for algebraic varieties defined over a base field
$k$ and taking values into an abelian, tensor category. The study of this
problem originated as the consequence of a general dissatisfaction
connected to an extensive use of topological methods in algebraic
geometry, with the result of producing several but insufficiently related
cohomological theories (Betti, de-Rham, \'etale etc.). The typical
example is furnished by a family of homomorphisms
$H^i_{\rm{et}}(X,\Q_\ell) \to H^i_{\rm{et}}(Y,\Q_\ell)$ connecting the
groups of \'etale cohomology of two (smooth, projective) varieties, as
the prime number $\ell$ varies, which are not connected, in general,
by any sort of (canonical) relation. 

\no The definition of a contravariant functor (functor of motivic cohomology)
\[
h: \mathcal V_k \to \mathcal M_k(\mathcal V_k),\qquad X\mapsto h(X)
\]
from the category $\mathcal V_k$ of
projective, smooth, irreducible algebraic varieties over $k$ to a semi-simple abelian category of
pure motives $\mathcal M_k(\mathcal V_k)$ is also tied up with the
definition of a universal cohomological theory through which every other
classical, cohomology $H^\cdot$ (here understood as contravariant functor) should factor by means of the introduction of
a fiber functor (realization $\otimes$-functor) $\omega$ connecting $\mathcal M_k(\mathcal V_k)$ to the abelian category
of (graded) vector-spaces over $\Q$
\begin{diagram}
\mathcal V_k &\rTo^{\qquad\quad H^\cdot\qquad\quad} &gr\mathcal V ect_\Q \\
&\rdTo_h &\uTo^\omega \\
& &\mathcal M_k(\mathcal V_k)
\end{diagram}
Following Grothendieck's
original viewpoint, the functor $h$
should implement the sought for mechanism of compatibilities (in \'etale
cohomology) and at the same time it should also describe a universal linearization of the theory of algebraic
varieties.

The definition of the category $\mathcal M_k(\mathcal V_k)$
arised from a classical construction in algebraic geometry which is based on the idea
of extending the collection of algebraic morphisms
in $\mathcal V_k$ by including the (algebraic) {\it correspondences}. A correspondence 
between two objects $X$ and $Y$ in $\mathcal V_k$ is
a multi-valued map which connects them. An algebraic correspondence is defined by means of an algebraic
cycle in the cartesian product $X\times Y$. The concept of
(algebraic) correspondence in geometry is much older than that of a
motive: it is in fact already present in several works of the
Italian school in algebraic geometry (\cf Severi's theory of
correspondences on algebraic curves). 

\no Grothendieck's new intuition was that the whole philosophy of motives is regulated by
the theory of  (algebraic) correspondences:

``{\it ...J'appelle motif sur $k$ quelque chose comme un groupe de
cohomologie $\ell$-adique d'un schema alg\'ebrique sur $k$, mais
consid\'er\'ee comme ind\'ependant de $\ell$, et avec sa structure
enti\`ere, ou disons pour l'instant sur $\Q$, d\'eduite de la
th\'eorie des cycles alg\'ebriques...}''
(\cf\cite{15}, Lettre 16.8.1964). \vspace{.05in}

Motives were envisioned with the hope to explain the intrinsic
relations between integrals of algebraic functions in one or more
complex variables. Their ultimate goal was to supply a machine
that would guarantee a generalization of the main results of
Galois theory to systems of polynomials equations in
several variables. Here, we refer in particular to
a higher-dimensional analog of the well-known result which describes the linearization of the
Galois-Grothendieck correspondence for the category $\mathcal V_k^o$ of \'etale, finite $k$-schemes
\[
\mathcal V_k^o \stackrel{\sim}{\to} \{\text{finite sets with \text{Gal}$(\bar k/k)$-action}\},\qquad X\mapsto X(\bar k)
 \]
by means of the equivalence between the
category of Artin motives and that of the representations of the absolute Galois group.

In the following section we shall describe how the
fundamental notions of the theory of motives arose from the study of
several classical problems in geometry and arithmetic.

\subsection{A first approach to motives}

A classical problem in algebraic geometry is that of computing the
solutions of a finite set of polynomial equations
\[
f_1(X_1,\ldots,X_m)=0,\ldots,f_r(X_1,\ldots,X_m)=0
\]
with coefficients in a finite field $\bF_q$. This study is naturally
formalized by introducing the generating series
\begin{equation}\label{zeta}
\zeta(X,t) = \exp(\sum_{m\ge 1}\nu_m\frac{t^m}{m})
\end{equation}
which is associated to the algebraic variety $X= V(f_1,\ldots,f_m)$ that is
defined as the set of the common zeroes of the polynomials
$f_1,\ldots,f_r$. 

\no Under the assumption that $X$ is smooth and
projective, the series \eqref{zeta} encodes the complete information
on the number of the rational points of the algebraic variety,
through the coefficients $\nu_m = |X(\bF_{q^m})|$. The integers
$\nu_m$ supply the cardinality of the set of the rational points of
$X$, computed in successive finite field extensions
$\bF_{q^m}$ of the base field $\bF_q$.

\no Intersection theory furnishes a general way to determine the number
$\nu_m$ as intersection number of two algebraic cycles on the
cartesian product $X\times X$: namely the diagonal $\Delta_X$ and
the graph $\Gamma_{Fr^m}$ of the $m$-th iterated composite of the
Frobenius morphism on the scheme $(X,\mathcal O_X)$:
\[
Fr: X\to X;\qquad Fr(P) = P,\quad f(\underline x) \mapsto f(\underline x^q),\quad\forall
f\in\mathcal O_X(U),~\forall~U\subset X~\text{open set}.
\]
The Frobenius endomorphism is in fact an interesting example of {\it
correspondence}, perhaps the most interesting one, for algebraic
varieties defined over finite fields. As a correspondence it induces a
commutative diagram
\begin{equation}\label{ac}
\begin{CD}
H^*(X\times X) @>{-\cap\Gamma_{Fr}}>> H^*(X\times X)\\
@AAp_1^*A @VV(p_2)_*V \\
H^*(X) @>Fr^*>> H^*(X)
\end{CD}
\end{equation}
in \'etale cohomology. Through the
commutativity of the above diagram one gets a way to express the action of the induced
homomorphism in cohomology, by means of the formula
\[
Fr^*(c) = (p_2)_*(p_1^*(c)\cap \Gamma_{Fr}),\qquad\forall c\in H^*(X),
\]
where $p_i: X\times X \to X$ denote the two projection maps. Here, `algebraic'
refers to the algebraic cycle $\Gamma_{Fr}\subset X\times
X$ that performs such a
correspondence. 

\no For particularly simple algebraic varieties, such as 
projective spaces $\bP^n$, the computation of the integers $\nu_m$
can be done by applying an elementary combinatorial argument
based on the set-theoretical description of the space $\bP^n(k) =
k^{n+1}\setminus\{0\}/k^\times$ ($k=$ any field). This has the effect
to produce the
interesting description
\begin{equation}\label{firstdec}
|\bP^n(\bF_{q^m})| = \frac{q^{m(n+1)}-1}{q^m-1} = 1 + q^m + q^{2m}
+\cdots+ q^{mn}.
\end{equation}
This decomposition of the set of the rational points of a
projective space was certainly a first source of inspiration in the
process of formalizing the foundations of the theory of motives. In
fact, one is naturally led to wonder on the casuality of the decomposition
\eqref{firstdec}, possibly ascribing such a result to the presence of a
cellular decomposition on the projective space which induces a related break-up on the set of the rational points. Remarkably, A.~Weil proved that a similar formula
holds also in the more general case of a smooth, projective
algebraic curve $C_{/{{\bF}_q}}$ of genus $g\ge 0$. In this case one shows that
\begin{equation}\label{deccurve}
|C(\bF_{q^m})| = 1 - \sum_{i=1}^{2g}\omega_i^m + q^m;\qquad
\omega_i\in\bar\Q,\quad |\omega_i| = q^{1/2}.
\end{equation}
These results  suggest that \eqref{firstdec} and
\eqref{deccurve} are the manifestation of a deep and intrinsic structure
that governs the geometry of  algebraic varieties. 

\no The development of
the theory of motives has in fact shown to us that this structure
reveals itself in several contexts: topologically, manifests its
presence in the decomposition of the cohomology $H^*(X) =
\oplus_{i\ge 0} H^i(X)$, whereas arithmetically it turns out that it is the same structure that controls the decomposition of the
series \eqref{zeta} as a rational function of $t$:
\[
\zeta(X,t) = \frac{\prod_{i\ge 0}
\det(1-tFr^*|H^{2i+1}_{et}(X))}{\prod_{i\ge
0}\det(1-tFr^*|H^{2i}_{et}(X))}.
\]
This is in fact a consequence of the description of the integers $\nu_m$ supplied by the
Lefschetz-Grothendieck trace formula (\cf~\cite{22})
\[
|X(\bF_{q^m})| = \sum_{i\ge 0}(-1)^i\text{tr}((Fr^m)^*|H^i_{et}(X)).
\]

\subsection{Grothendieck's pure motives}

Originally, Grothendieck proposed a general framework for
a so called category of {\it numerically effective motives} ${\bf M}(k)_\Q$
over a field $k$ and with rational coefficients. This
category is defined by enlarging  the  category $\mathcal V_k$ of
smooth, projective algebraic varieties over $k$ (and algebraic
morphisms) by following the so-called procedure of pseudo-abelian envelope.
This construction is performed in two steps: at first one enlarges
the set of morphisms of $\mathcal V_k$ by including 
(rational) algebraic correspondences of degree zero, modulo numerical equivalence, then one
performs a pseudo-abelian envelope by formally including among the
objects, kernels of idempotent morphisms.

\no Let us assume for simplicity that the algebraic varieties are
irreducible (the general case is then deduced from this, by additivity).
For any given $X,Y\in\text{Obj}(\mathcal V_k)$, one works with correspondences $f: X
\dashrightarrow Y$ which are  elements of codimension
equal to $\dim X$  in the rational graded algebra
\[A^*(X\times Y) = C^*(X\times Y)\otimes\Q/\sim_{num}
\]
of algebraic cycles modulo {\it numerical equivalence}. We recall
that two algebraic cycles on an algebraic variety $X$ are said to be
numerically equivalent $Z\sim_{num} W$, if
\begin{equation}\label{degree}
\text{deg}(Z\cdot T) = \text{deg}(W\cdot T),
\end{equation}
 for any algebraic cycle
$T$ on $X$. Here, by $\text{deg}(V)$ we mean the degree of the
algebraic cycle $V = \sum_{\text{finite}} m_\alpha V_\alpha\in
C^*(X)$.

\no The degree defines a homomorphism from the free abelian group of
algebraic cycles $C^*(X) = \oplus_i C^i(X)$ to the integers. On the
components $C^i(X)$, the map is defined as follows
\[
\text{deg}: C^i(X) \to \Z,\qquad \text{deg}(V) =
\begin{cases} \sum_\alpha m_\alpha& \text{if $i=\dim X$},\\0& \text{if $i<\dim
X$}.\end{cases}
\]
The symbol `$~\cdot~$' in \eqref{degree} refers to the intersection
product structure on $C^*(X)$, which is well-defined under the assumption
of proper intersection. If $Z\cap T$ is proper (\ie
$\text{codim}(Z\cap T) = \text{codim}(Z) + \text{codim}(T)$), then
intersection theory supplies the definition of an intersection
cycle $Z\cdot T\in C^*(X)$. Moreover, the intersection product is commutative
and associative whenever is defined.

\no Passing from the free abelian group $C^*(X)$ to the quotient
$C^*(X)/\sim$, modulo a suitable equivalence relation on cycles,
allows one to use classical results of algebraic geometry (so called
Moving Lemmas) which lead to the definition of a ring structure.
One then defines intersection cycle classes in
general, even when cycles do not intersect properly, by intersecting
equivalent cycles which fulfill the required geometric property of proper intersection.

\no It is natural to guess that the use of the numerical equivalence in the original definition of the
category of motives was motivated by the study of classical constructions in enumerative geometry,
such as for example the computation of the number of the
rational points of an algebraic variety defined over a finite field. One of the main 
original goals was 
to show that for a suitable definition of an equivalence relation on
algebraic cycles, the corresponding category of motives is {\it
semi-simple}. This means that the objects $M$ in the
category decompose, following the rules of a {\it theory of
weights} (\cf section~\ref{fundstructures}), into direct {\it factors} $M =
\oplus_i M_i(X)$, with $M_i(X)$ simple (\ie indecomposable)
motives associated to smooth, projective algebraic varieties. The
importance of achieving such a result is quite evident if one seeks, for
example, to understand categorically the decomposition $H^*(X) =
\oplus_i H^i(X)$ in cohomology, or if one wants to recognize the role of motives
in the factorization of zeta-functions of algebraic
varieties.\smallskip

\no Grothendieck concentrated his efforts on the numerical
equivalence relation which is the coarsest among the equivalence
relations on algebraic cycles. So doing, he attacked the problem of
the semi-simplicity of the category of motives from the easiest side.
However, despite a promising departing point, the statement on the
semi-simplicity escaped all his efforts. In fact, the result he was able to
reach at that time was dependent on the assumption of the {\it
Standard Conjectures}, two strong topological statements on
algebraic cycles. The proof of the semi-simplicity of the
category of motives for numerical equivalence (the only equivalence
relation producing this result) was achieved only much later on in
the development of the theory (\cf\cite{24}). The proof
found by U. Jannsen uses a fairly elementary but ingenious idea
which mysteriously eluded Grothendieck's intuition as well as all the
mental grasps of several mathematicians after him.

By looking at the construction of (pure) motives in perspective, one
immediately recognizes the predominant role played by the morphisms
over the objects, in the category ${\bf M}(k)_\Q$. This was
certainly a great intuition of Grothendieck. This idea led  to a
systematic study of the properties of algebraic cycles and their
decomposition by means of {\it algebraic projectors}, that is
algebraic cycles classes $p\in A^{\dim X}(X\times X)$ satisfying the
property
\[p^2 =
p\circ p = p.
\]
Notice that in order to make sense of the notion of a projector and
more in general, in order to define a  law of composition
`$\circ$' on algebraic correspondences, one needs to use the
ring structure on the graded algebra $A^*$. The operation `$\circ$' is
defined as follows. Let us assume for simplicity, that the
algebraic varieties are connected (the general case can be easily
deduced from this). Then, two algebraic correspondences $f_1\in
A^{\dim X_1+i}(X_1\times X_2)$ (of degree $i$) and $f_2\in A^{\dim
X_2+j}(X_2\times X_3)$ (of degree $j$) compose accordingly to the
following rule (bi-linear, associative)
\[
A^{\dim X_1+i}(X_1\times X_2)\times A^{\dim X_2+j}(X_2\times X_3)
\to A^{\dim X_1+i+j}(X_1\times X_3)
\]
\[
(f_1,f_2) \mapsto f_2\circ f_1 = (p_{13})_*((p_{12}^*(f_1)\cdot
(p_{23})^*(f_2))).
\]
In the particular case of projectors $p: X \dashrightarrow X$, one
is restricted, in order to make a sense of the condition
$p\circ p = p$, to use only particular types of algebraic
correspondences: namely those of degree zero. These are the elements of
the abelian group $A^{\dim X}(X\times X)$.

The objects of the category ${\bf M}(k)_\Q$ are then pairs $(X,p)$,
with $X\in\text{Obj}(\mathcal V_k)$ and $p$ a projector. This way, one
attains  the notion of a $\Q$-linear, pseudo-abelian,
monoidal category (the $\otimes$-monoidal structure is deduced from
the cartesian product of algebraic varieties), together with the
definition of a contravariant functor
\[
h: \mathcal V_k \to {\bf M}(k)_\Q,\qquad X\mapsto h(X) =
(X,\text{id}).
\]
Here $(X,\text{id})$ denotes the motive associated to $X$ and
$\text{id}$ means the trivial (\ie identity) projector
associated to the diagonal $\Delta_X$. More in general, $(X,p)$
refers to the motive $ph(X)$ that is cut-off on $h(X)$ by 
the (range of the) projector $p: X \dashrightarrow X$. Notice that
images of projectors are formally included among the objects of
${\bf M}(k)_\Q$, by the procedure of the pseudo-abelian envelope.
The cut-off performed by a projector $p$ on the space determines a corresponding
operation in cohomology (for any classical Weil
theory), by singling out the sub-vector space $pH^*(X)\subset
H^*(X)$.

The category ${\bf M}(k)_\Q$ has two important basic objects: ${\bf
1}$ and ${\bf L}$. ${\bf 1}$ is the {\it unit motive}
\[
{\bf 1} = (\Spec(k),\text{id}) = h(\Spec(k)).
\]
This is defined by the zero-dimensional algebraic variety
associated to a point, whereas
\[
{\bf L} = ({\bf P}^1,\pi_2), \quad \pi_2 = {\bf P}^1\times
\{P\},\quad P\in {\bf P}^1(k)
\]
is the so-called {\it Lefschetz motive}. This motive determines, jointly
with ${\bf 1}$, a decomposition of the motive associated to the
projective line ${\bf P}^1$
\begin{equation}\label{decp1}
h({\bf P}^1) = {\bf 1}\oplus {\bf L}.
\end{equation}
One can show that the algebraic cycles ${\bf P}^1\times\{P\}$ and $\{P\}\times{\bf
P}^1$ on ${\bf P}^1\times{\bf P}^1$ do not depend on the choice of the
rational point $P\in{\bf P}^1(k)$ and that their sum is a cycle
equivalent to the diagonal. This fact implies that the decomposition
\eqref{decp1} is canonical. More in general, it follows from the K\"unneth decomposition
of the diagonal $\Delta$ in ${\bf P}^n\times {\bf P}^n$ by algebraic cycles 
$\Delta = \pi_0 +\cdots + \pi_n$, (\cf~\cite{31} and \cite{18} for the
details) that the motive of a projective space ${\bf
P}^n$ decomposes into pieces (simple motives)
\begin{equation}\label{dec}
h({\bf P}^n) = h^0({\bf P}^n)\oplus h^2({\bf P}^n)\cdots\oplus
h^{2n}({\bf P}^n)
\end{equation}
where $h^{2i}({\bf P}^n) = ({\bf P}^n, \pi_{2i}) = (h^2({\bf
P}^n))^{\otimes i}$, $\forall i>0$. It is precisely this
decomposition which implies the decomposition
\eqref{firstdec} on the rational points, when $k =
{\bF}_q$!

\no For (irreducible) curves, and in the presence of a
rational point $x\in C(k)$, one obtains  a similar decomposition (non
canonical)
\[h(C) = h^0(C)\oplus h^1(C)\oplus h^2(C)
\]
with $h^0(C) = (C,\pi_0=\{x\}\times C)$, $h^2(C) =
(C,\pi_2=C\times\{x\})$ and $h^1(C) = (C,1-\pi_0-\pi_2)$. This
decomposition is responsible for the formula \eqref{deccurve}.

In fact, one can prove that these decompositions partially generalize to any object
$X\in\text{Obj}(\mathcal V_k)$. In the presence of a rational point,
or more in general by choosing a positive zero-cycle $Z =
\sum_\alpha m_\alpha Z_\alpha\in C^{\dim X}(X\times X)$ (here $X$ is
assumed irreducible for simplicity and $\dim X = d$), one constructs
two rational algebraic cycles
\[
\pi_0 = \frac{1}{m}(Z\times X),\quad \pi_{2d} = \frac{1}{m}(X\times
Z);\qquad m = \deg(Z) = \sum_\alpha m_\alpha>0
\]
which determine two projectors $\pi_0,\pi_{2d}$ in the Chow group
$CH^d(X\times X)\otimes\Q$ of rational algebraic cycles modulo
rational equivalence. The corresponding classes in $A^d(X\times X)$
(if not zero) determine two motives $(X,\pi_0)\simeq h^0(X)$ and
$(X,\pi_{2d})\simeq h^{2d}(X)$ (\cf\eg\cite{34}).

For the applications, it is convenient to enlarge the category of
effective motives by formally adding the tensor product inverse
${\bf L}^{-1}$ of the Lefschetz motive: one usually refers to it
as to the {\it Tate motive}. It corresponds, from the more
refined point of view of Galois theory, to the cyclotomic
characters. This enlargement of ${\bf M}(k)_\Q$ by the so-called
``virtual motives'' produces an abelian, semi-simple category
$\mathcal M_k(\mathcal V_k)_\Q$ of pure motives for numerical
equivalence. The objects of this category are now triples $(X,p,m)$,
with $m\in\Z$. Effective motives are of course objects of this
category and they are described by triples $(X,p,0)$. The Lefschetz
motive gains a new interpretation in this category as ${\bf L} =
(\Spec(k),\text{id},-1)$. The Tate motive is defined by ${\bf
L}^{-1} = (\Spec(k),\text{id},1)$ and is therefore reminissent of
(in fact induces) the notion of Tate structure $\Q(1)$ in Hodge
theory.

In the category $\mathcal M_k(\mathcal V_k)_\Q$, the set of
morphisms connecting two motives $(X,p,m)$ and $(Y,q,n)$
 is defined by
\[
\text{Hom}((X,p,m),(Y,q,n)) = q\circ A^{\dim X-m+n}(X\times Y)\circ
p.
\]
In particular, $\forall f=f^2\in \text{End}((X,p,m))$, one defines
the two motives 
\[
\Image(f) = (X,p\circ f\circ p,m),\qquad \Ker(f) =
(X,p-f,m).
\]
These determine a canonical decomposition of any virtual motive as
$\Image(f)\oplus \Image(1-f) \stackrel{\sim}{\to} (X,p,m)$, where
the direct sum of two motives as the above ones is defined by
the formula
\[
(X,p,m)\oplus (Y,q,m) = (X\coprod Y,p+q,m).
\]
The general definition of the direct sum of two motives requires a  bit more of
formalism which escapes this short overview: we refer to \op for the
details. 

\no The tensor structure
\[
(X,p,m)\otimes(Y,q,n) = (X\times Y,p\otimes q,m+n)
\]
and the involution (\ie auto-duality) which is defined, for $X$ 
irreducible and $\dim X = d$ by the functor
\[
\vee: (\mathcal M_k(\mathcal V_k)_\Q)^{op} \to \mathcal M_k(\mathcal
V_k)_\Q,\qquad (X,p,m)^\vee = (X,p^t,d-m)
\]
(the general case follows from this by
applying additivity), determine the structure of a {\it rigid} $\otimes$-{\it category} on $\mathcal
M_k(\mathcal V_k)_\Q$. Here, $p^t$ denotes the transpose
correspondence associated to $p$ (\ie the transpose of the graph).
One finds, for example, that ${\bf L}^\vee = {\bf L}^{-1}$. In the
particular case of the effective motive $h(X)$, with $X$ irreducible
and $\dim X = d$, this involution determines the notion of Poincar\'e
duality
\[
h(X)^\vee = h(X)\otimes {\bf L}^{\otimes(-d)}
\]
that is an auto-duality which induces the Poincar\'e duality isomorphism in any
classical cohomological theory.

\subsection{Fundamental structures}\label{fundstructures}

A category of pure motives over a field $k$ and with coefficients in
a field $K$ (of characteristic zero) is supposed to satisfy, to be
satisfactory, several basic properties and to be endowed with a few
fundamental structures. In the previous section we have described
the historically-first example of a category of pure motives  and we
have reviewed some of its basic properties (in that case $K = \Q$).
One naturally wonders about the description of others categories of motives associated
to finer (than the numerical) equivalence relations on algebraic
cycles: namely the categories of motives for homological or rational or algebraic
equivalence relations. However, if
one seeks to work with a semi-simple category, the afore mentioned result of Jannsen tells
us that the numerical equivalence is the only adequate relation. The semi-simplicity
property is also attaint if one assumes
Grothendieck's Standard Conjectures. Following the report
of Grothendieck in \cite{20}, these conjectures arose from the hope to
understand (read prove) the conjectures of Weil on the zeta-function
of an algebraic variety defined over a finite field. It was well-known to
Grothendieck that the Standard Conjectures imply the Weil's
Conjectures. These latter statements became a theorem in the early
seventies (1974), only a few years later the time when
Grothendieck stated  the Standard Conjectures
(1968-69). The proof by Deligne of the Weil's conjectures, however,
does not make any use of the Standard Conjectures, these latter questions remain
still unanswered at the present time. The moral lesson seems to be that
geometric topology and the theory of algebraic
cycles govern in many central aspects the foundations of algebraic
geometry.

The Standard Conjectures of ``Lefschetz type'' and of ``Hodge type''
are stated in terms of algebraic cycles modulo homological
equivalence (\cf\cite{26}). They imply two further important
conjectures. One of these states the equality of the homological and
the numerical equivalence relations, the other one, of
``K\"unneth type'' claims that the K\"unneth decomposition of the
diagonal in cohomology, can be described by means of (rational)
algebraic cycles. It is generally accepted, nowadays, to refer to the full set of
the four conjectures, when one quotes the Standard Conjectures.
In view of their expected consequences, one is naturally led to
study a category of pure {\it motives for homological equivalence}. In
fact, there are several candidates for this category since the
definition depends upon the choice of a Weil cohomological theory
(\ie Betti, \'etale, de-Rham, crystalline, etc) with coefficients in
a field $K$ of characteristic zero.

Let us fix a cohomological theory $X\mapsto H^*(X) = H^*(X,K)$ for
algebraic varieties in the category $\mathcal V_k$. Then, the
construction of the corresponding category of motives $\mathcal
M_k^{hom}(\mathcal V_k)_K$ for homological equivalence is given following
a procedure similar to that we have explained earlier on in this paper
for the category of motives for numerical equivalence
and with rational coefficients. The only difference is that now morphisms in the
category $\mathcal M_k^{hom}(\mathcal V_k)_K$ are defined by means
of algebraic correspondences modulo homological equivalence. At this point,
one makes explicit use of the axiom ``cycle map'' that characterizes (together
with finiteness,
Poincar\'e duality, K\"unneth formula, cycle map, weak and strong
Lefschetz theorems) any Weil cohomological theory (\cf\cite{26}).

\no The set of algebraic morphisms connecting objects in the category
$\mathcal V_k$ is enlarged by including  multi-valued maps
$X\dashrightarrow Y$ that are defined as a $K$-linear combination
of elements of the vector spaces
\[
C^*(X\times Y)\otimes_\Z F/\sim_{hom},
\]
where $F\subset K$ is a subfield. Two cycles $Z,W\in C^*(X\times
Y)\otimes F$ are {\it homologically equivalent} $Z\sim_{hom} W$ if
their image, by means of the cycle class map
\[
\gamma: C^*(X\times Y)\otimes F \to H^*(X\times Y)
\]
is the same. This leads naturally to the definition of a
subvector-space $A^*_{hom}(X\times Y)\subset H^*(X\times Y)$
generated by the image of the cycle class map $\gamma$. These spaces
define the correspondences in the category $\mathcal
M_k^{hom}(\mathcal V_k)_K$. If $X$ is purely $d$-dimensional, then
\[
\text{Corr}^r(X,Y) := A^{d+r}_{hom}(X\times Y).
\]
In general, if $X$ decomposes into several connected components 
$X = \coprod_i X_i$,  one
sets $\text{Corr}^r(X,Y) = \oplus_i\text{Corr}^r(X_i,Y)$. In direct
analogy to the construction of correspondences for numerical
equivalence, the ring structure (``cap-product'') in cohomology
determines a composition law `$\circ$' among correspondences.

The category $\mathcal M_k^{hom}(\mathcal V_k)_K$ is then defined as
follows: the objects are triples $M=(X,p,m)$, where
$X\in\text{Obj}(\mathcal V_k)$, $m\in\Z$ and $p=p^2\in
\text{Corr}^0(X,X)$ is an idempotent. The collection of morphisms
between two motives $M=(X,p,m)$, $N=(Y,q,n)$ is given by the set
\[
\text{Hom}(M,N) = q\circ\text{Corr}^{n-m}(X,Y)\circ p.
\]
This procedure determines a pseudo-abelian, $K$-linear tensor
category. The tensor law is given by the formula
\[
(X,p,m)\otimes (Y,q,n) = (X\times Y,p\times q,m+n).
\]
The commutativity and associativity constraints are induced by the
obvious isomorphisms $X\times Y \stackrel{\sim}{\to} Y\times X$,
$X\times(Y\times Z)\stackrel{\sim}{\to}(X\times Y)\times Z$. The
unit object in the category is given by ${\bf 1} =
(\text{Spec}(k),\text{id},0)$. One shows that $\mathcal M_k^{hom}(\mathcal V_k)_K$
is a {\it rigid} category, as it is endowed with an auto-duality
functor
\[
\vee: \mathcal M_k^{hom}(\mathcal V_k)_K \to (\mathcal
M_k^{hom}(\mathcal V_k)_K)^{op}.
\]
For any object $M$, the functor $-\otimes M^\vee$ is
left-adjoint to $-\otimes M$ and $M^\vee\otimes -$ is right-adjoint
to $M\otimes -$. In the case of an irreducible
variety $X$, the internal Hom  is defined by the motive
\[
\underline{\text{Hom}}((X,p,m),(Y,q,n)) = (X\times Y,p^t\times
q,\dim(X)-m+n).
\]

The Standard Conjecture of K\"unneth type (which is assumed from now
on in this section) implies that the K\"unneth components of the diagonal
$\pi^i_X\in\text{Corr}^0(X,X)$ determine a complete system of
orthogonal, central idempotents. This important statement implies
that  the motive $h(X)\in \mathcal M_k^{hom}(\mathcal V_k)_K$ has
the expected direct sum decomposition (unique)
\[
h(X) = \bigoplus_{i=0}^{2\dim X}h^i(X),\qquad h^i(X) = \pi^i_Xh(X).
\]
The cohomology functor $X\mapsto H^*(X)$ factors through the
projection $h(X)\to h^i(X)$. More in general, one shows that every motive $M=
(X,p,m)$ gets this way a $\Z$-{\it grading} structure by setting
\begin{equation}\label{grading}
(X,p,m)^r = (X,p\circ\pi^{r+2m},r).
\end{equation}
This grading is respected by all morphisms in the category and
defines the structure of a {\it graduation by weights} on the
objects. On a motive $M = (X,p,m)= ph(X)\otimes{\bf L}^{\otimes(-m)}
= ph(X)(m)$ in the category, one sets

\[
M = \oplus_i\text{Gr}_i^w(M),\qquad \text{Gr}_i^w(M) =
ph^{2m+i}(X)(m)
\]
where $\text{Gr}_i^w(M)$ is a pure motive of weight $i$. One finds
for example, that ${\bf 1}$ has weight zero, ${\bf
L}=(\text{Spec}(k),1,-1)$ has weight $2$ and that ${\bf L}^{-1} =
(\text{Spec}(k),1,1)$ has weight $-2$. More in general, the motive
$M=(X,p,m) = ph(X)(m)$ has weight $-2m$.

In order to achieve further important properties, one needs to
modify the natural commutativity constraint
$\psi=\oplus_{r,s}\psi^{r,s}$, $\psi^{r,s}: M^r\otimes N^s
\stackrel{\sim}{\to} N^s\otimes M^r$, by defining
\begin{equation}\label{mod}
\psi_{new}=\oplus_{r,s}(-1)^{rs}\psi^{r,s}.
\end{equation}

\no We shall denote by $\tilde{\mathcal M}_k^{hom}(\mathcal V_k)_K$ the
category of motives for homological equivalence in which one has
implemented the modification \eqref{mod}  on the tensor product
structure.

An important structure on a category of pure motives (for homological
equivalence) is given by
assigning to an object $X\in\text{Obj}(\mathcal M_k^{hom}(\mathcal V_k)_K)$ 
a {\it motivic cohomology} $H^i_{mot}(X)$.
$H^i_{mot}(X)$ is a pure motive of weight $i$. This way, one
views pure motives as a universal cohomological theory for algebraic
varieties. The main property of the motivic cohomology is
that it defines a universal realization of any  given Weil cohomology
theory $H^*$. Candidates for these motivic
cohomology theories have been proposed by A. Beilinson \cite{3},
in terms of eigenspaces of Adams operations in algebraic $K$-theory
\ie $H^{2j-n}_{mot}(X,\Q(j)) = K_n(X)^{(j)}$ and by S. Bloch
\cite{5}, in terms of higher Chow groups \ie
$H^{2j-n}_{mot}(X,\Q(j)) = CH^j(X,n)\otimes\Q$.

The assignment of a Weil cohomological theory with coefficients in a
field $K$ which contains an assigned  field $F$ is equivalent to the definition
of an exact {\it realization $\otimes$-functor} of $\tilde{\mathcal
M}_k^{hom}(\mathcal V_k)_K$ in the category of $K$-vector spaces
\begin{equation}\label{realiz}
r_{H^*}: \tilde{\mathcal M}_k^{hom}(\mathcal V_k)_K \to
\text{Vect}_K,\qquad r_{H^*}(H^i_{mot}(X)) \simeq H^i(X).
\end{equation}
In particular, one obtains the realization $r_{H^*}({\bf L}^{-1}) = H^2({\bf P}^1)$
which defines the notion of the {\it Tate twist} in cohomology. More precisely:\vspace{.05in}

\no -- in \'etale cohomology: $H^2({\bf P}^1) = \Q_\ell(-1)$, where
$\Q_\ell(1):= \displaystyle{\varprojlim_m} ~\mu_{\ell^m}$ is a
$\Q_\ell$-vector space of dimension one endowed with the cyclotomic
action of the absolute Galois group $G_k = \text{Gal}(\bar k/k)$.
The ``twist'' (or torsion) $(r)$ in \'etale cohomology corresponds
to the torsion in Galois theory defined by the $r$-th power of the
cyclotomic character (Tate twist)\smallskip

\no -- in de-Rham theory: $H^2_{DR}({\bf P}^1) = k$, with the Hodge
filtration defined by $F^{\le 0}=0$, $F^{>0} = k$. Here, the effect
of the torsion $(r)$ is that of shifting the Hodge filtration of
$-r$-steps (to the right)\smallskip

\no -- in Betti theory: $H^2({\bf P}^1) = \Q(-1) := (2\pi i)^{-1}\Q$.
The bi-graduation on $H^2({\bf P}^1)\otimes\C \simeq \C$ is purely
of type $(1,1)$. The torsion $(r)$ is here identified with the
composite of a homothety given by a multiplication by $(2\pi
i)^{-r}$ followed by a shifting by $(-r,-r)$ of the Hodge
bi-graduation.

Using the structure of rigid tensor-category one introduces
the notion of {\it rank} associated to a motive $M= (X,p,m)$
in $\tilde{\mathcal M}_k^{hom}(\mathcal V_k)_K$. The rank of $M$ is defined as
the trace of $\text{id}_M$ \ie the trace of the morphism
$\epsilon\circ\psi_{new}\circ\eta\in\text{End}({\bf 1})$, where
\[
\epsilon: M\otimes M^\vee \to {\bf 1},\qquad \eta: {\bf 1} \to
M^\vee\otimes M
\]
are {\it resp.} the evaluation and co-evaluation morphisms
satisfying $\epsilon\otimes\text{id}_M\circ\text{id}_M\otimes\eta =
\text{id}_M$,
$\text{id}_{M^\vee}\otimes\epsilon\circ\eta\otimes\text{id}_{M^\vee}
= \text{id}_{M^\vee}$. In general, one sets
\[
\text{rk}(X,p,m) = \sum_{i\ge 0}\dim pH^i(X)\ge 0.
\]
Under the assumption of the Standard Conjectures (more precisely
under the assumption that homological and numerical equivalence
relations coincide) and that $\text{End}({\bf 1}) = F$
($\text{char}(F)=0$), the tannakian formalism invented by
Grothendieck and developed by Saavedra \cite{33}, and Deligne
\cite{16} implies that the abelian, rigid, semi-simple tensor
category $\tilde{\mathcal M}_k^{hom}(\mathcal V_k)_K$ is endowed
with an exact, faithful $\otimes$-fibre functor to the category of
graded $K$-vector spaces
\begin{equation}\label{ff}
\omega: \tilde{\mathcal M}_k^{hom}(\mathcal V_k)_K \to
\text{VectGr}_K,\qquad \omega(H^*_{mot}(X)) = H^*(X)
\end{equation}
which is compatible with the realization functor. This formalism
defines a {\it tannakian} ({\it neutral} if $K=F$) structure on the
category of motives. One then introduces the {\it tannakian group}
\[
G = \underline{\text{Aut}}^{\otimes}(\omega)
\]
as a $K$-scheme in affine groups.  Through the tannakian
formalism one shows that the fibre functor $\omega$
realizes an equivalence of rigid tensor categories
\[
\omega:\tilde{\mathcal M}_k^{hom}(\mathcal V_k)_K
\stackrel{\sim}{\to} \text{Rep}_F(G),
\]
where $\text{Rep}_F(G)$ denotes the rigid tensor category of finite
dimensional, $F$ representations of the tannakian group $G$. 
This way, one establishes a quite useful dictionary between
categorical $\otimes$-properties and properties of the associated
groups. Because we have assumed all along the Standard Conjectures, the
semi-simplicity of the category $\tilde{\mathcal M}_k^{hom}(\mathcal
V_k)_K$ implies that $G$ is an algebraic, pro-reductive group, that is
$G$ is the projective limit of reductive $F$-algebraic
groups.

\no The tannakian theory is a linear analog of the theory of  finite,
\'etale coverings of a given connected scheme. This theory was developed by
Grothendieck in SGA1 (theory of the pro-finite $\pi_1$). For this
reason the group $G$ is usually referred to as the {\it motivic
Galois group} associated to $\mathcal V_k$ and $H^*$. In the case of
algebraic varieties of dimension zero (\ie for Artin motives) the
tannakian group $G$ is nothing but the (absolute) Galois
group $\text{Gal}(\bar k/k)$.

In any reasonable cohomological theory the functors $X\mapsto
H^*(X)$ are deduced by applying standard methods of homological
algebra to the related derived functors $X \mapsto R\Gamma(X)$
which associate to an object in $\mathcal V_k$ a bounded complex of
$k$-vector spaces, in a suitable triangulated category $\mathcal
D(k)$ of complexes of modules over $k$, whose heart is the category
of motives. This is the definition of cohomology as
\[
H^i(X) = H^iR\Gamma(X).
\]
Under the assumption that the functors $R\Gamma$ are  realizations
of corresponding motivic functors \ie $R\Gamma =
r_{H^*}R\Gamma_{mot}$, one expects the existence of a
(non-canonical) isomorphism in $\mathcal D(k)$
\begin{equation}\label{motcoh}
R\Gamma_{mot}(X) \simeq \oplus_i H^i_{mot}(X)[-i].
\end{equation}
Moreover, the introduction of the motivic derived functors
$R\Gamma_{mot}$ suggests the definition of the following groups of
{\it absolute cohomology}
\[
H^i_{abs}(X) = \text{Hom}_{\mathcal D(k)}({\bf
1},R\Gamma_{mot}(X)[i]).
\]
For a general motive $M = (X,p,m)$, one defines
\begin{equation}\label{ExtM}
H^i_{abs}(M) =  \text{Hom}_{\mathcal D(k)}({\bf 1},M[i]) =
\text{Ext}^i({\bf 1},M).
\end{equation}
The motives $H^i_{mot}(X)$ and the groups of absolute motivic
cohomology are related by a spectral sequence
\[
E_2^{p,q} = H^p_{abs}(H^q_{mot}(X)) ~\Rightarrow~ H^{p+q}_{abs}(X).
\]

\subsection{Examples of pure motives}

The first interesting examples of pure motives arise by considering
the category $\mathcal V_k^o$ of \'etale, finite $k$-schemes. An
object in this category is a scheme $X=\text{Spec}(k')$, where $k'$
is a commutative $k$-algebra of finite dimension which satisfies the
following properties. Let $\bar k$ denote a fixed separable closure
of $k$
\begin{enumerate}
\item[1.]~$k'\otimes\bar k \simeq \bar k^{[k':k]}$
\item[2.]~$k' \simeq \prod k_\alpha$, for $k_\alpha/k$ finite,
separable field extensions
\item[3.]~$|X(\bar k)|=[k':k]$.
\end{enumerate}
The corresponding rigid, tensor-category of motives with
coefficients in a field $K$ is usually referred to as the category
of {\it Artin motives}: $\mathcal C\mathcal V^o(k)_K$. 

\no The definition of this category is independent of the choice of the
equivalence relation on cycles as the objects of $\mathcal V_k^o$
are smooth, projective $k$-varieties of dimension zero. One also
sees that passing from $\mathcal V_k^o$ to $\mathcal C\mathcal
V^o(k)_K$ requires adding new objects in order to attain the
property that the category of motives is abelian. One can verify
this already for $k = K= \Q$, by considering the real quadratic
extension $k' = \Q(\sqrt 2)$ and the one-dimensional non-trivial
representation of $G_\Q = \text{Gal}(\bar\Q/\Q)$ that factors
through the character of order two of $\text{Gal}(k'/\Q)$. This
representation does not correspond to any object in $\mathcal
V_\Q^o$, but can be obtained as the image of the projector
$p=\frac{1}{2}(1-\sigma)$, where $\sigma$ is the generator of
$\text{Gal}(k'/\Q)$. Therefore,
$\text{image}(p)\in\text{Obj}(\mathcal C\mathcal V^o(\Q)_\Q)$ is a
new object.

The category of Artin motives is a semi-simple, $K$-linear, monoidal
$\otimes$-category. When $\text{char}(K)=0$, the commutative diagram
of functors
\begin{equation*}
\begin{CD}
\mathcal V_k^o @>GG>\sim> \{\text{sets with $\text{Gal}(\bar
k/k)$-continuous action}\}\\
@VVhV @VVlV\\
\mathcal C\mathcal V^o(k)_K @>(*)>\sim> \{\text{finite dim.
$K$-v.spaces with linear $\text{Gal}(\bar k/k)$-continuous
action}\}
\end{CD}
\end{equation*}
where $l$ is the contravariant functor of linearization
\[
S\mapsto K^S,\qquad (g(f))(s) = f(g^{-1}(s)),\quad\forall
g\in\text{Gal}(\bar k/k)
\]
determines a linearization of the Galois-Grothendieck correspondence
(GG) by means of the equivalence of categories $(*)$. This is
provided by the fiber functor
\[
\omega: X \to H^0(X_{\bar k},K) = K^{X(\bar k)}
\]
and by applying the tannakian formalism. It follows that $\mathcal
C\mathcal V^o(k)_K$ is $\otimes$-equivalent to the category
$\text{Rep}_K\text{Gal}(\bar k/k)$ of representations of the
absolute Galois group $G_k = \text{Gal}(\bar k/k)$.

\no These results were the departing point for Grothendieck's
speculations on the definition of higher dimensional Galois
theories (\ie Galois theories associated to system of polynomials in
several variables) and for the definition of the corresponding
motivic Galois groups.

\section{Endomotives: an overview}

The notion of an {\it endomotive} in
noncommutative geometry (\cf\cite{7}) is the natural generalization
of the classical concept of an Artin motive for the noncommutative spaces 
which are defined by semigroup actions on projective limits of 
zero-dimensional algebraic varieties, endowed with an action of the 
absolute Galois group. This notion applies quite naturally for instance,
to the study of several examples of quantum statistical dynamical systems
whose time evolution describes important number-theoretic properties
of a given field $k$ (\cf{\it op.cit}, \cite{14}).  

\no There are two distinct definitions of an endomotive: one speaks of
algebraic or analytic endomotives depending upon the context and the
applications. 

\no When $k$ is a number field, there is a functor
connecting the two related categories. Moreover, the abelian category of
Artin motives embeds naturally as a full subcategory in the category
of algebraic endomotives (\cf~Theorem~\ref{NCArtinThm}) and this
result motivates the statement that the theory of endomotives defines a
natural generalization of the classical theory of (zero-dimensional)
 Artin motives.

In noncommutative geometry, where the properties of a space
(frequently highly singular from a classical viewpoint) are analyzed
in terms of the properties of the associated noncommutative algebra and its (space of)
irreducible representations, it is quite natural to look for a
suitable abelian category which enlarges the original, non-additive
category of algebras and in which one may also apply the standard
techniques of homological algebra. Likewise in the construction of
a theory of motives, one seeks to work within a triangulated
category endowed with several structures. These include for
instance, the definition of (noncommutative) motivic objects playing
the role of motivic cohomology (\cf~\eqref{motcoh}), the
construction of a universal (co)homological theory representing in
this context the absolute motivic cohomology (\cf~\eqref{ExtM}) and
possibly also the set-up of a noncommutative tannakian formalism to 
motivate in rigorous mathematical terms the presence of certain universal groups
of symmetries associated to renormalizable quantum field theories
(\cf\eg~\cite{12}).

\no A way to attack these problems is that of enlarging the original
category of algebras and morphisms by introducing a ``derived''
category of modules enriched with a suitable notion of
correspondences connecting the objects that should also account for the
structure of Morita equivalence which represents the noncommutative
generalization of the notion of isomorphism for commutative
algebras.

\subsection{The abelian category of cyclic
modules}\label{cyclicmodules}

The sought for enlargement of the category $Alg_k$ of (unital)
$k$-algebras and (unital) algebra homomorphisms is defined by
introducing a new category $\Lambda_k$ of cyclic $k(\Lambda)$-modules.
The objects of this category are modules over the {\it cyclic category}
$\Lambda$. This latter has the same objects as the simplicial
category $\Delta$ ($\Lambda$ contains $\Delta$ as sub-category). We
recall that an object in $\Delta$ is a totally ordered set
\[
[n] = \{0<1<\ldots<n\}
\]
for each $n\in\N$, and a morphism
\[
f: [n] \to [m]
\]
is described by an order-preserving map of sets $f: \{0,1,\ldots,n\}
\to \{0,1,\ldots,m\}$. The set of morphisms in $\Delta$ is generated
by faces $\delta_i: [n-1] \to [n]$ (the injection that misses $i$)
and degeneracies $\sigma_j: [n+1] \to [n]$ (the surjection which
identifies $j$ with $j+1$) which satisfy several standard simplicial
identities (\cf\eg\cite{9}). The set of morphisms in $\Lambda$
is enriched by introducing a new collection of morphisms: the {\it
cyclic morphisms}. For each $n\in\N$, one sets
\[
\tau_n: [n] \to [n]
\]
fulfilling the relations
\begin{equation}\label{ae}
\begin{array}{lll} \tau_n  \delta_i = \delta_{i-1} \tau_{n-1} & 1 \leq i \leq
n , & \tau_n  \delta_0 = \delta_n \\[3mm] \tau_n  \sigma_i =
\sigma_{i-1} \tau_{n+1} &1 \leq i \leq n , &\tau_n  \sigma_0 =
\sigma_n
\tau_{n+1}^2 \\[3mm] \tau_n^{n+1} = 1_n. &&  \end{array}
\end{equation}
The objects of the category $\Lambda_k$ are $k$-modules over
$\Lambda$ (\ie $k(\Lambda)$-modules). In categorical language this
means functors
\[
\Lambda^{op} \to Mod_k
\]
($Mod_k=$ category of $k$-modules). Morphisms of
$k(\Lambda)$-modules are therefore natural transformations between
the corresponding functors.

\no It is evident that $\Lambda_k$ is an abelian category, because
of the  interpretation of a morphism in $\Lambda_k$ as a collection
of $k$-linear maps of $k$-modules $A_n \to B_n$
($A_n,B_n\in\text{Obj}(Mod_k)$) compatible with faces, degeneracies
and cyclic operators. Kernels and cokernels of these morphisms
define objects of the category $\Lambda_k$, since their definition
is given point-wise.

To an algebra $\mathcal A$ over a field $k$, one associates the
$k(\Lambda)$-module $\mathcal A^\natural$. For each $n\ge 0$ one
sets
\[
\cA^{\natural}_n = \underbrace{\cA\otimes \cA\otimes\cdots \otimes
\cA}_{\text{(n+1)-times}}.
\]
The {\it cyclic morphisms} on $\cA^{\natural}$ correspond to the
cyclic permutations of the tensors, while the face and the degeneracy
maps correspond to the algebra product of consecutive tensors and
the insertion of the unit. This construction determines a functor
\[
\natural: Alg_k \to \Lambda_k
\]
{\it Traces} $\varphi: \cA \to k$ give rise naturally to  morphisms
\[
\varphi^\natural:\cA^{\natural}\to k^\natural,\qquad
\varphi^\natural(a_0\otimes \cdots \otimes a_n)=\,\varphi(a_0
\cdots a_n)
\] 
in $\Lambda_k$.
The main result of this construction is the following canonical
description of the cyclic cohomology of an algebra $\mathcal A$ over
a field $k$ as the derived functor of the functor which assigns to a
$k(\Lambda)$-module its space of traces
\begin{equation}\label{ExtC}
HC^n(\cA) = \Ext^n (\cA^\natural, k^\natural)
\end{equation}
(\cf\cite{9},\cite{29}). This formula is the analog of
\eqref{ExtM}, that describes the absolute motivic cohomology group
of a classical motive as an Ext-group computed in a triangulated
category of motives $\mathcal D(k)$. In the present context, on the other hand, the
derived groups ${\rm Ext}^n$ are taken in the abelian category of
$\Lambda_k$-modules.

\no The description of the cyclic cohomology as a derived functor in the
cyclic category determines a useful procedure to embed the nonadditive category
of algebras and algebra homomorphisms in the ``derived'' abelian
category of $k(\Lambda)$-modules. This construction provides a natural framework for
the definition of the objects of a category of noncommutative
motives.

\no Likewise in the construction of the category of motives, one is faced
with the problem of finding the ``motivated maps'' connecting
cyclic modules. The natural strategy is that of enlarging
the collection of cyclic morphisms which are functorially induced by
homomorphisms between (noncommutative) algebras, by implementing an
adequate definition of (noncommutative) correspondences. The notion
of an algebraic correspondence in algebraic geometry, as a multi-valued
map defined by an algebraic
cycle modulo a suitable equivalence relation, has here an analog
with the notion of Kasparov's bimodule and the associated class in
$KK$-theory (\cf \cite{25}). Likewise in classical motive theory,
one may prefer to work with (compare) several versions of
correspondences. One may decide to retain the full information supplied by a
group action on a given algebra (\ie a noncommutative space) rather
than partially loosing this information by moding out with the equivalence
relation (homotopy in $KK$-theory).

\subsection{Bimodules and $KK$-theory}

There is a natural way to associate a cyclic morphism to a (virtual)
correspondence and hence to a class in
$KK$-theory. Starting with the category of separable $C^*$-algebras
and $*$-homomorphisms, one enlarges the collection of morphisms
connecting two unital algebras $\mathcal A$ and $\mathcal B$, by
including correspondences defined by elements of Kasparov's
bivariant $K$-theory
\[
\Hom(\cA,\cB)=KK(\cA,\cB)
\]
(\cite{25}, \cf also \S 8 and \S 9.22 of \cite{4}). More
precisely, correspondences are defined by Kasparov's bimodules,
that means by triples
\[\mathcal E = \mathcal E(\mathcal A,\mathcal B)=(E,\phi,F)
\]
which satisfy the following conditions:\vspace{.05in}

\no --~$E$ is a countably generated Hilbert module over $\cB$

\no --~$\phi$ is a $*$-homomorphism of $\cA$ to bounded linear operators
on $E$ (\ie $\phi$ gives $E$ the structure of an $\cA$-$\cB$ bimodule)

\no --~$F$ is a bounded linear operator on $E$ such that the operators
$[F,\phi(a)]$, $(F^2-1)\phi(a)$, and $(F^*-F)\phi(a)$ are compact
for all $a\in \cA$.

\no A Hilbert module $E$ over $\cB$ is a right $\cB$-module with a {\it
positive}, $\cB$-valued inner product which satisfies $\langle x, yb
\rangle= \langle x, y\rangle b$, $\forall x,y\in E$ and $\forall
b\in \cB$, and with respect to which $E$ is complete (\ie complete
in the norm $\parallel x\parallel = \sqrt{\parallel\langle
x,x\rangle\parallel}$).

\no Notice that Kasparov bimodules are Morita-type of correspondences.
They generalize $*$-homomorphisms of $C^*$-algebras since the latter ones
may be re-interpreted as Kasparov bimodules of the form
$(\cB,\phi,0)$.

\no Given a Kasparov's bimodule $\mathcal E = \mathcal E(\mathcal
A,\mathcal B)$, that is a $\mathcal A$-$\mathcal B$ Hilbert bimodule
$E$ as defined above, one associates, under the assumption that
$E$ is a projective $\mathcal B$-module of finite type
(\cf\cite{7}, Lemma~2.1), a cyclic morphism
\[
\mathcal E^\natural\in \text{Hom}(\mathcal A^\natural,\mathcal
B^\natural).
\]
This result allows one to define an enlargement of the
collection of cyclic morphisms in the category $\Lambda_k$ of
$k(\Lambda)$-modules, by considering Kasparov's projective bimodules
of finite type, as correspondences.

\no One then implements the {\it homotopy equivalence relation} on the
collection of Kasparov's bimodules. Two Kasparov's modules are said to be
homotopy equivalent $(E_0,\phi_0,F_0)\sim_h (E_1,\phi_1,F_1)$ if
there is an element
\[
(E,\phi,F) \in \mathcal E(\cA,I\cB),\quad I\cB=\{ f:[0,1]\to
\cB\,|\, f \text{ continuous} \}
\]
which performs a unitary homotopy deformation between the two
modules. This means that $(E\hat\otimes_{f_i}\cB,f_i\circ
\phi,f_i(F))$ is unitarily equivalent to $(E_i,\phi_i,F_i)$ or
equivalently re-phrased, that there is a unitary in bounded
operators from $E\hat\otimes_{f_i}\cB$ to $E_i$ intertwining the
morphisms $f_i\circ \phi$ and $\phi_i$ and the operators $f_i(F)$
and $F_i$. Here $f_i: I\cB\to \cB$ is the evaluation at the
endpoints.

There is a binary operation on the set of all Kasparov $\mathcal
A$-$\mathcal B$ bimodules, given by the direct sum. By definition, the
group of Kasparov's bivariant $K$-theory is the set of homotopy
equivalence classes $c(\mathcal E(\mathcal A,\mathcal B))\in
KK(\mathcal A,\mathcal B)$ of Kasparov's modules $\mathcal
E(\cA,\cB)$. This set has a natural structure of abelian group with
addition induced by direct sum.

This bivariant version of $K$-theory is reacher than both $K$-theory
and $K$-homology, as it carries an intersection product. There is a
natural  bi-linear, associative composition (intersection) product
\[
\otimes_{\mathcal B}: KK(\mathcal A,\mathcal B) \times KK(\mathcal
B,\mathcal C) \to KK(\mathcal A,\mathcal C)
\]
for all $\cA, \cB$ and $\cC$ separable $C^*$-algebras. This product
is compatible with composition of morphisms of $C^*$-algebras.

\no $KK$-theory is also endowed with a bi-linear, associative exterior
product
\[
\otimes: KK(\cA,\cB)\otimes KK(\cC,\cD)\to KK(\cA\otimes
\cC,\cB\otimes \cD),
\]
which is defined in terms on the composition product by
\[
c_1\otimes c_2 = (c_1\otimes
1_{\cC})\otimes_{\cB\otimes\C}(c_2\otimes 1_{\cB}).
\]
A slightly different formulation of $KK$-theory, which simplifies
the definition of this external tensor product is obtained by
replacing in the data $(E,\phi,F)$ the operator $F$ by an {\em
unbounded}, regular self-adjoint operator $D$. The corresponding $F$
is then given by $D(1+D^2)^{-1/2}$ (\cf \cite{1}).

The above construction which produces an enlargement of the category
of separable $C^*$-algebras by introducing correspondences as
morphisms determines an additive, although non abelian category
$\cK\cK$ (\cf \cite{4} \S 9.22.1). This category is also known to
have a triangulated structure (\cf\cite{32}) and
this result is in agreement with the construction of the triangulated
category $\mathcal D(k)$ in motives theory, whose heart is expected
to be the category of (mixed) motives (\cf
section~\ref{fundstructures} and \cite{17}). A more refined
analysis based on the analogy with the construction of a category of
motives suggests that one should probably perform a further enlargement by
passing to the pseudo-abelian envelope of $\cK\cK$, that is by
formally including among the objects also ranges of idempotents in
$KK$-theory.

In section~\ref{anendo} we will review the category of analytic
endomotives where maps are given in terms of \'etale
correspondences described by  spaces $\mathcal Z$ arising
from locally compact \'etale groupoids $\mathcal G = \mathcal
G(X_\alpha,S,\mu)$ associated to zero-dimensional, singular quotient
spaces $X(\bar k)/S$ with associated $C^*$-algebras $C^*(\cG)$. In
view of what we have said in this section, it would  be also possible
to define a category where morphisms are given by classes $c(\cZ)\in
KK(C^*(\cG), C^*(\cG'))$ which describe sets of equivalent triples
$(E,\phi,F)$, where $(E,\phi)$ is given in terms of a bimodule $\cM_\cZ$
with the trivial grading $\gamma=1$ and the zero endomorphism $F=0$.
The definition of the category of
analytic endomotives is more refined because the definition of the
maps in this category does not require to divide by homotopy
equivalence.

The comparison between correspondences for motives given by
algebraic cycles and correspondences for noncommutative spaces given
by bimodules (or elements in $KK$-theory) is particularly easy in
the zero-dimensional case because the equivalence relations play no
role. Of course, it would be quite interesting to investigate the higher
dimensional cases, in view of a unified framework for
motives and noncommutative spaces which is suggested for example, by
the recent results on the Lefschetz trace formula for archimedean local
factors of $L$-functions of motives (\cf~\cite{7}, Section~7).

\no A way to attack this problem is by comparing the notion of a
correspondence given by an  algebraic cycle with the notion
of a geometric correspondence used in topology (\cf\cite{2},
\cite{13}). For example, it is easy to see that the definition of
an algebraic correspondence can be reformulated as a particular case
of the topological (geometrical) correspondence and it is also known
that one may associate to the latter a class in $KK$-theory. In the
following two sections, we shall review and comment on these ideas.

\subsection{Geometric correspondences}

In geometric topology, given two smooth manifolds $X$ and $Y$ (it is enough 
to assume that  $X$ is a locally compact parameter
space), a {\it topological (geometric) correspondence} is given by
the datum
\[
X\stackrel{f_X}{\leftarrow} (Z,E) \stackrel{g_Y}{\to} Y
\]
where:\vspace{.05in}

\no --~$Z$ is a smooth manifold

\no --~$E$ is a complex vector bundle over $Z$

\no --~$f_X:Z\to X$ and $g_Y:Z\to Y$ are continuous maps, with $f_X$
proper and $g_Y$ $K$-oriented (orientation in $K$-homology).

\no Unlike in the definition of an algebraic correspondence
(\cf~section~\ref{fundstructures}) one does not require that $Z$ is
a subset of the cartesian product $X\times Y$. This  flexibility is
balanced by the implementation of the extra piece of datum given by
the vector bundle $E$.
To any such correspondence $(Z,E,f_X,g_Y)$ one associates a class
in Kasparov's K-theory
\begin{equation}\label{kZEfg} c(Z,E,f_X,g_Y) = (f_X)_* (
(E)\otimes_Z
(g_Y)! )\in KK(X,Y).
\end{equation}
$(E)$ denotes the class of $E$ in $KK(Z,Z)$ and  $(g_Y)!$ is
the element in $KK$-theory which fulfills the Grothendieck
Riemann-Roch formula. 

\no We recall that given two smooth manifolds
$X_1$ and $X_2$ and a continuous oriented map $f: X_1 \to X_2$, the
element $f! \in KK(X_1,X_2)$ determines the Grothendieck
Riemann--Roch formula
\begin{equation}\label{GRR}
{\rm ch}(F \otimes f!) =f_! ({\rm Td}(f) \cup {\rm ch}(F)),
\end{equation}
for all $F\in K^*(X_1)$, with ${\rm Td}(f)$ the Todd genus
\begin{equation}\label{Toddf}
{\rm Td}(f) = {\rm Td}(TX_1)/{\rm Td}(f^* TX_2).
\end{equation}

\no The composition of two correspondences $(Z_1,E_1,f_X,g_Y)$ and
$(Z_2,E_2,f_Y,g_W)$ is given  by taking the fibered product
$Z=Z_1\times_Y Z_2$ and the bundle $E=\pi_1^*E_1 \times \pi_2^*
E_2$, with $\pi_i:Z\to Z_i$ the projections. This determines the
composite correspondence $(Z,E,f_X,g_W)$. In fact, one also needs to
assume a {\it transversality} condition on the maps $g_Y$ and $f_Y$
in order to ensure that the fibered product $Z$ is a smooth
manifold. The homotopy invariance of both $g_Y!$ and $(f_X)_*$ show
however that the assumption of transversality is `generically'
satisfied.

\no Theorem 3.2 of \cite{13} shows that  Kasparov product in
$KK$-theory `$\otimes$' agrees with the composition of
correspondences, namely
\begin{equation}\label{prodGeomKK}
c(Z_1,E_1,f_X,g_Y)\otimes_Y c(Z_2,E_2,f_Y,g_W) =c(Z,E,f_X,g_W) \ \ \
\in KK(X,W).
\end{equation}

\subsection{Algebraic correspondences and $K$-theory}

In algebraic geometry, the notion of correspondence that comes
closest to the definition of a geometric correspondence (as an
element in $KK$-theory) is obtained by considering classes of
algebraic cycles in algebraic $K$-theory (\cf \cite{30}).

Given two smooth and projective algebraic varieties $X$ and $Y$, we
denote by $p_X$ and $p_Y$ the projections of $X\times Y$ onto $X$
and $Y$ respectively and we assume that they are {\it proper}. Let
$Z\in C^*(X\times Y)$ be an algebraic cycle. For simplicity, we
shall assume  that $Z$ is irreducible (the general case follows by
linearity). We denote by $f_X=p_X|_Z$ and $g_Y=p_Y|_Z$ the
restrictions of  $p_X$ and $p_Y$ to $Z$.

\no To the irreducible subvariety $T\stackrel{i}{\hookrightarrow} Y$ one
naturally associates the coherent $\cO_{Y}$-module $i_* \cO_T$. For
simplicity of notation we write it as $\cO_T$. We use a similar
notation for the coherent sheaf $\mathcal O_Z$, associated to the
irreducible subvariety $Z\hookrightarrow X\times Y$. Then, the sheaf
pullback
\begin{equation}\label{XYmod}
p_{Y}^* \cO_T =p_{Y}^{-1} \cO_T \otimes_{p_{Y}^{-1}\cO_{Y}}
\cO_{X\times Y}
\end{equation}
has a natural structure of  $\cO_{X\times Y}$-module. The map on
sheaves that corresponds to the cap product by $Z$ on cocycles is
given by
\begin{equation}\label{cyclemap2}
Z: \cO_T \mapsto (p_X)_* \left( p_Y^* \cO_T \otimes_{\cO_{X\times
Y}} \cO_Z  \right).
\end{equation}
Since $p_X$ is proper, the resulting sheaf is coherent. Using
\eqref{XYmod}, we can write equivalently
\begin{equation}\label{cyclemap3}
Z: \cO_T \mapsto (p_X)_* \left( p_Y^{-1} \cO_T \otimes_{p_Y^{-1}
\cO_Y} \cO_Z  \right).
\end{equation}

\no We recall that the functor $f_!$ is the right adjoint to $f^*$ (\ie
$f^*f_!=id$) and that $f_!$ satisfies the Grothendieck Riemann--Roch
formula
\begin{equation}\label{GRR2}
{\rm ch}(f_!(F)) =f_! ({\rm Td}(f) \cup {\rm ch}(F)).
\end{equation}
Using this result, we can equally compute the intersection product of
\eqref{cyclemap2} by first computing
\begin{equation}\label{shriekZ}
\cO_T \otimes_{\cO_Y} (p_Y)_! \cO_Z
\end{equation}
and then applying $p_Y^*$. Using \eqref{GRR2} and \eqref{GRR} we
know that we can replace \eqref{shriekZ} by $\cO_T \otimes_{\cO_Y}
(\cO_Z \otimes (p_Y)!)$ with the same effect in $K$-theory.

\no Thus, to a correspondence in the sense of \eqref{cyclemap2} that is defined
by the image in $K$-theory of an algebraic cycle $Z\in C^*(X\times
Y)$ we associate the geometric class
\[
\mathcal F(Z) = c(Z,E,f_X,g_Y) \in KK(X,Y),
\]
with $f_X=p_X|_Z$, $g_Y=p_Y|Z$ and with the bundle
$E=\cO_Z$.

\no The composition of correspondences is given in terms of the
intersection product of the associated cycles. Given three smooth
projective varieties $X$, $Y$ and $W$ and (virtual) correspondences
$U=\sum a_i Z_i\in C^*(X\times Y)$ and $V=\sum c_j Z_j '\in
C^*(Y\times W)$, with $Z_i\subset X\times Y$ and $Z_j '\subset
Y\times W$ closed reduced irreducible subschemes, one defines
\begin{equation}\label{intprod}
U \circ V= (\pi_{13})_* ((\pi_{12})^* U \bullet (\pi_{23})^* V).
\end{equation}
$\pi_{12}:X\times Y\times W \to X\times Y$, $\pi_{23}:
X\times Y\times W \to Y\times W$, and $\pi_{13}: X\times Y\times W
\to X\times W$ denote, as usual, the projection maps.

Under the assumption of `general position' which is the algebraic analog of
the transversality requirement in topology, we obtain the following
result

\begin{proposition}[\cf\cite{7} Proposition~6.1]\label{corrcorr2}
Suppose given three smooth projective varieties $X$, $Y$, and $W$ and 
algebraic correspondences $U$ given by  $Z_1\subset X\times Y$ and $V$
given by  $Z_2\subset Y\times W$. Assume that $(\pi_{12})^*
Z_1$ and $(\pi_{23})^* Z_2$ are in general position in $X\times
Y\times W$. Then assigning to a cycle $Z$ the topological
correspondence $\cF(Z)=(Z,E,f_X,g_Y)$ satisfies
\begin{equation}\label{prodprod}
\cF(Z_1\circ Z_2)=\cF(Z_1)\otimes_Y \cF(Z_2),
\end{equation}
where $Z_1\circ Z_2$ is the product of algebraic cycles and
$\cF(Z_1)\otimes_Y \cF(Z_2)$ is the Kasparov product of the
topological correspondences.
\end{proposition}

\no Notice that, while in the topological (smooth) setting
transversality can always be achieved by a small deformation (\cf \S
III, \cite{13}), in the algebro-geometric framework one needs to
modify the above construction if the cycles
are not in general position. In this case the formula
\[
[\cO_{T_1}]\otimes [\cO_{T_2}]= [\cO_{T_1\circ T_2}]
\]
which describes the product in $K$-theory in terms of the
intersection product of algebraic cycles must be modified by
implementing  Tor-classes and one works with a product defined by
the formula (\cite{30}, Theorem 2.7)
\begin{equation}\label{prodKint2}
[\cO_{T_1}]\otimes [\cO_{T_2}]= \sum_{i=0}^n (-1)^i \left[{\rm
Tor}_i^{\cO_X} (\cO_{T_1},\cO_{T_2})\right].
\end{equation}

\section{Algebraic endomotives}

To define the category of {\it algebraic endomotives} one replaces
the category $\mathcal V^o_k$ of reduced, finite-dimensional
commutative algebras (and algebras homomorphisms) over a field $k$
by the category of noncommutative algebras (and algebras
homomorphisms) of the form
\[
\mathcal A_k = A\rtimes S.
\]
$A$ denotes a unital algebra which is an inductive limit of
commutative algebras $A_\alpha\in\text{Obj}(\mathcal V^o_k)$. $S$ is
a unital, {\it abelian semigroup} of algebra endomorphisms
\[\rho: A \to A.
\]
Moreover, one imposes the condition that for $\rho\in S$, $e=\rho(1)\in A$ is an
{\it idempotent} of the algebra and that $\rho$ is an isomorphism of
$A$ with the compressed algebra $eAe$.

\no The crossed product algebra $\mathcal A_k$ is defined by formally
adjoining to $A$ new generators $U_\rho$ and $U_\rho^*$, for
$\rho\in S$, satisfying the algebraic rules
\begin{alignat}{3}\label{algrulesS}
U_\rho^*U_\rho &= 1, \quad & U_\rho U_\rho^* &=
\rho(1),\qquad\forall\rho\in S,\\\notag U_{\rho_1\rho_2} &=
U_{\rho_1}U_{\rho_2}, \quad & U_{\rho_2\rho_1}^*
&= U_{\rho_1}^*U_{\rho_2}^*,\quad\forall\rho_j\in S,\\
U_\rho x &= \rho(x)U_\rho, \quad &  xU_{\rho}^* &=
U_\rho^*\rho(x),\quad\forall\rho\in S,~\forall x\in A.\notag
\end{alignat}
Since $S$ is abelian, these rules suffice to show that $\mathcal
A_k$ is the linear span of the monomials $U_{\rho_1}^*aU_{\rho_2}$,
for $a\in A$ and $\rho_j\in S$.

\no Because $A = \displaystyle{\varinjlim_\alpha} A_\alpha$, with
$A_\alpha$ reduced, finite-dimensional commutative algebras over
$k$, the construction of $\mathcal A_k$ is in fact determined by
assigning a {\it projective system} $\{X_\alpha\}_{\alpha\in I}$ of
varieties in $\mathcal V^o_k$ ($I$ is a countable indexing set),
with $\xi_{\beta,\alpha}: X_\beta\to X_\alpha$ morphisms in
$\mathcal V^o_k$ and with a suitably defined action of $S$. Here, we
have implicitly used the equivalence between the category of finite
dimensional commutative $k$-algebras and the category of affine
algebraic varieties over $k$.

\no The graphs $\Gamma_{\xi_{\beta,\alpha}}$ of the connecting morphisms
of the system  define $G_k = \text{Gal}(\bar k/k)$-invariant subsets
of $X_\beta(\bar k)\times X_\alpha(\bar k)$ which in turn describe
$\xi_{\beta,\alpha}$ as algebraic correspondences. We denote by
\[
X = \varprojlim_\alpha X_\alpha,\qquad \xi_\alpha: X \to X_\alpha
\]
the associated {\it pro-variety}. The compressed algebra $eAe$
associated to the idempotent $e = \rho(1)$ determines a subvariety
$X^e\subset X$ which is in fact isomorphic to $X$, via the induced
morphism $\tilde\rho: X \to X^e$.

\no The noncommutative space defined by $\mathcal A_k$ is the quotient
of $X(\bar k)$ by the action of $S$, \ie of the action of the
$\tilde\rho$'s.

\no The Galois group $G_k$ acts on $X(\bar k)$ by composition. By
identifying  the elements of $X(\bar k)$ with characters, \ie with
$k$-algebra homomorphisms $\chi: A \to \bar k$, we write the action
of $G_k$ on $A$ as
\begin{equation}\label{charaction}
\alpha(\chi) = \alpha\circ\chi: A \to \bar k,\qquad\forall\alpha\in
G_k,~\forall\chi\in X(\bar k).
\end{equation}
This action commutes with the maps $\tilde\rho$, \ie
$(\alpha\circ\chi)\circ\rho = \alpha\circ(\chi\circ\rho)$. Thus the
whole construction of the system $(X_\alpha,S)$ is
$G_k$-equivariant. This fact does not mean however, that $G_k$ acts
by automorphisms on $\mathcal A_k$!

\no Moreover, notice  that the algebraic construction of the crossed-product
algebra $\mathcal A_k$ endowed with the actions of $G_k$ and $S$ on
$X(\bar k)$ makes sense also when $\text{char}(k)>0$.

\no When $\text{char}(k)=0$, one defines the set of correspondences
$M(\mathcal A_k,\mathcal B_k)$ by using the notion of Kasparov's
bimodules $\mathcal E(\mathcal A_k,\mathcal B_k)$ which are
projective and finite  as right modules. This way, one obtains a
first realization of the resulting category of noncommutative
zero-dimensional motives in the abelian category of
$k(\Lambda)$-modules.

In general, given  $(X_\alpha,S)$, with $\{X_\alpha\}_{\alpha\in I}$
a projective system of Artin motives and $S$ a semigroup of
endomorphisms of $X=\displaystyle{\varprojlim_\alpha}X_\alpha$ as
above,  the datum  of the semigroup action is encoded naturally by
the {\it algebraic groupoid}
\[
\mathcal G = X\rtimes S.
\]
This is defined in the following way. One considers the Grothendieck
group $\tilde S$ of the abelian semigroup $S$. By using the
injectivity of the partial action of $S$, one may also assume that
$S$ embeds in $\tilde S$. Then, the action of $S$ on $X$ extends to
define a {\it partial action} of $\tilde S$. More precisely, for $s
= \rho_2^{-1}\rho_1\in\tilde S$ the two projections
\[
E(s) = \rho_1^{-1}(\rho_2(1)\rho_1(1)),\qquad F(s) =
\rho_2^{-1}(\rho_2(1)\rho_1(1))
\]
only depend on $s$ and the map $s: A_{E(s)} \to A_{F(s)}$ defines an
isomorphism of reduced algebras. It is immediate to verify that
$E(s^{-1}) = F(s) = s(E(s))$ and that $F(ss') \ge F(s)s(F(s'))$. The
algebraic groupoid $\mathcal G$ is defined as the disjoint union
\[
\mathcal G = \coprod_{s\in\tilde S}X^{F(s)}
\]
which corresponds to the commutative direct-sum of reduced algebras
\[
\bigoplus_{s\in\tilde S}A_{F(s)}.
\]
The range and the source maps in $\mathcal G$ are given \resp by the
natural projection from $\mathcal G$ to $X$ and by its composition
with the antipode ${\bf S}$ which is defined, at the algebra level,
by ${\bf S}(a)_s = s(a_{s^{-1}}),~\forall s\in\tilde S$. The
composition in the groupoid corresponds to the product of monomials
$aU_sbU_t = as(b)U_{st}$.

Given two systems $(X_\alpha,S)$ and $(X_{\alpha '}',S')$, with
associated crossed-product algebras $\mathcal A_k$ and $\mathcal
B_k$ and groupoids $\mathcal G = \mathcal G(X_\alpha,S)$ and
$\mathcal G' = \mathcal G(X'_{\alpha'},S')$  a {\it geometric
correspondence} is given by a $(\mathcal G,\mathcal G')$-space
$Z=\text{Spec}(C)$, endowed with a right action of $\mathcal G'$
which fulfills the following {\it \'etale} condition. Given a space
such as $\mathcal G'$, that is a disjoint union  of zero-dimensional
pro-varieties over $k$, a right action of $\mathcal G'$ on $Z$ is
given by a map $g: Z \to X'$ and a collection of partial
isomorphisms
\begin{equation}\label{groupoidact11}
z\in g^{-1}(F(s)) \mapsto z\cdot s\in g^{-1}(E(s))
\end{equation}
fulfilling the following rules for partial action of the abelian
group $\tilde S$
\begin{equation}\label{groupoidact12}
g(z\cdot s) = g(z)\cdot s,\qquad z\cdot(ss') = (z\cdot s)\cdot
s'\quad\text{on}~g^{-1}(F(s)\cap s(F(s'))).
\end{equation}
Here $x\mapsto x\cdot s$ denotes the partial action of $\tilde S$ on
$X'$. One checks that such an action gives to the $k$-linear space
$C$ a structure of right module over $\mathcal B_k$. The action of
$\mathcal G'$ on $Z$ is \'etale if the corresponding module $C$ is
{\it finite and projective} over $\mathcal B_k$.

Given two systems $(X_\alpha,S)$ and $(X'_{\alpha '},S')$ as above,
an {\it \'etale correspondence} is therefore a $(\mathcal
G(X_\alpha,S),\mathcal G(X_{\alpha'}',S'))$-space $Z$ such that the
right action of $\mathcal G(X_{\alpha'}',S')$ is \'etale.

The $\Q$-linear space of {\it (virtual) correspondences}
\[
\text{Corr}((X_\alpha,S),(X_{\alpha'}',S'))
\]
is the rational vector space of formal linear combinations $U =
\sum_ia_iZ_i$ of \'etale correspondences $Z_i$, modulo the relations
arising from isomorphisms and equivalences: $Z\coprod Z' \sim Z +
Z'$. The composition of correspondences is given by the fiber
product over a groupoid. Namely, for three systems $(X_\alpha,S)$,
$(X'_{\alpha'},S')$, $(X''_{\alpha''},S'')$ joined by
correspondences
\[
(X_\alpha,S) \leftarrow Z\to (X'_{\alpha'},S'),\qquad
(X'_{\alpha'},S')\leftarrow W \to (X''_{\alpha''},S''),
\]
their composition is given by the rule
\begin{equation}\label{comprule}
Z\circ W = Z\times_{\mathcal G'}W
\end{equation}
that is the fiber product over the groupoid $\mathcal G' = \mathcal
G(X'_{\alpha'},S')$.

Finally, a system $(X_\alpha,S)$ as above is said to be {\it
uniform} if the normalized counting measures $\mu_\alpha$ on
$X_\alpha$ satisfy $\xi_{\beta,\alpha}\mu_\alpha = \mu_\beta.$

\begin{definition} The category $\mathcal E\mathcal V^o(k)_K$ of {\it
algebraic endomotives} with coefficients in a fixed extension $K$ of
$\Q$ is the (pseudo)abelian category generated by the following
objects and morphisms. The objects are uniform systems
$M=(X_\alpha,S)$ of Artin motives over $k$, as above. The set of
morphisms in the category connecting two objects $M = (X_\alpha,S)$
and $M' = (X_{\alpha'}',S')$ is defined as
\[
\text{Hom}(M,M') =
\text{Corr}((X_\alpha,S),(X_{\alpha'}',S'))\otimes_\Q K.
\]
\end{definition}

The category $\mathcal C\mathcal V^o(k)_K$ of Artin motives embeds
as a full sub-category in the category of algebraic endomotives
\[
\iota: \mathcal C\mathcal V^o(k)_K \to \mathcal E\mathcal V^o(k)_K.
\]
The functor $\iota$ maps an Artin motive $M = X$ to the system
$(X_\alpha,S)$ with $X_\alpha = X$, $\forall\alpha$ and $S =
\{\text{id}\}$.

\subsection{Examples of algebraic endomotives}\label{genendo}

The category of algebraic endomotives is inclusive of a large and general
class of examples of noncommutative spaces $\mathcal A_k = A\rtimes
S$ which are described by semigroup actions on projective systems of
Artin motives.

\no One may consider, for instance a pointed algebraic variety $(Y,y_0)$
over a field $k$ and a countable, unital, abelian semigroup $S$ of
{\it finite} endomorphisms of $(Y,y_0)$, {\it unramified} over
$y_0\in Y$. Then, there is a  system $(X_s,S)$ of Artin motives
over $k$ which is constructed from these data. More precisely, for $s\in S$, one sets
\begin{equation}\label{Xs}
X_s=\{ y\in Y:\, s(y)=y_0 \}.
\end{equation}
For a pair $s,s'\in S$, with $s'=sr$, the connecting map
$\xi_{s,s'}: X_{sr}\to X_s$ is defined by
\begin{equation}\label{Xsr}
X_{sr} \ni y \mapsto r(y)\in X_s.
\end{equation}
This is an example of a system indexed by the semigroup $S$ itself,
with partial order given by divisibility. One sets
$X=\displaystyle{\varprojlim_s} X_s$.

\no Since $s(y_0)=y_0$, the base point $y_0$ defines a component $Z_s$
of $X_s$ for all $s\in S$. The pre-image $\xi_{s,s'}^{-1}(Z_s)$  in
$X_{s'}$ is a union of components of $X_{s'}$. This defines a
projection $e_s$ onto an open and closed subset $X^{e_s}$ of the
projective limit $X$.

\no It is easy to see that the semigroup $S$ acts on the projective
limit $X$ by partial isomorphisms $\beta_s: X \to X^{e_s}$ defined
by the property
\begin{equation}\label{endoS}
\beta_s: X \to X^{e_s},\quad \xi_{su}(\beta_s(x))=\xi_u(x),\qquad
\forall u\in S, \forall x\in X.
\end{equation}

The map $\beta_s$ is well-defined since the set $\{ su:\, u\in S \}$
is cofinal and $\xi_u(x)\in X_{su}$, with $su \xi_u(x)=s(y_0)=y_0$.
The image of $\beta_s$ is in $X^{e_s}$, since by definition of
$\beta_s$: $\xi_s(\beta_s(x))=\xi_1(x)=y_0$. For $x\in X^{e_s}$, we
have $\xi_{su}(x)\in X_u$. This shows that $\beta_s$ defines an
isomorphism of $X$ with $X^{e_s}$, whose inverse map is given by
\begin{equation}\label{invrhos}
\xi_u(\beta_s^{-1}(x))=\xi_{su}(x), \, \forall x\in X^{e_s}, \forall
u\in S.
\end{equation}

The corresponding algebra morphisms $\rho_s$ are then given by
\begin{equation}\label{endoS1}
\rho_s(f)(x)=f(\beta_s^{-1}(x)),\, \forall x\in X^{e_s}, \ \
\rho_s(f)(x)=0 ,\, \forall x\notin X^{e_s}\,.
\end{equation}

This class of examples  also fulfill an {\it equidistribution
property}, making the uniform normalized counting measures $\mu_s$
on $X_s$ compatible with the projective system and inducing a
probability measure on the limit $X$. Namely, one has
\begin{equation}\label{measlim}
\xi_{s',s} \mu_s = \mu_{s'},\, \forall s,s'\in S.
\end{equation}

\section{Analytic endomotives}\label{anendo}

In this section we assume that $k$ is a number field. We fix an
embedding $\sigma: k\hookrightarrow \C$ and we denote by $\bar k$ 
an algebraic closure of $\sigma(k)\subset \C$ in $\C$.

\no When taking points over $\bar k$, algebraic endomotives yield
$0$-dimensional singular quotient spaces $X(\bar k)/S$, which can be
described by means of locally compact \'etale groupoids $\cG(\bar k)$
and the associated crossed product $C^*$-algebras $C(X(\bar
k))\rtimes S$. This construction gives rise to the category of {\it
analytic endomotives}.

One starts off by considering a uniform system $(A_\alpha,S)$ of
Artin motives over $k$ and the algebras
\begin{equation}\label{complexify}
A_\C= A \otimes_k \C= \varinjlim_\alpha A_\alpha\otimes_k \C ,
 \    \    \   \cA_\C= \cA_k\otimes_k \C= A_\C \rtimes S .
\end{equation}
The assignment
\begin{equation}\label{embed}
a \in A \to \hat a, \ \ \ \hat a(\chi)=\,\chi(a)\ \ \ \forall \chi
\in X = \varprojlim_\alpha X_\alpha
\end{equation}
defines an involutive embedding of algebras $A_\C\subset C(X)$. The
$C^*$-completion $C(X)$ of $A_\C$ is an abelian AF $C^*$-algebra.
One sets
\[
\bar\cA_\C=C(X)\rtimes S.
\]
This is the $C^*$-completion of the algebraic crossed product
$A_\C\rtimes S$. It is defined by the algebraic relations
\eqref{algrulesS} with the involution which is apparent in the
formulae (\cf\cite{27},\cite{28}).

In the applications that require  to work with cyclic (co)homology,
it is important to be able to restrict from $C^*$-algebras such as
$\bar\cA_\C$ to canonical dense subalgebras
\begin{equation}\label{densesub}
\cA_\C= C^\infty(X)\rtimes_{alg} S \subset \bar\cA_\C
\end{equation}
where $C^\infty(X)\subset C(X)$ is the algebra of locally constant
functions. It is to this category of smooth algebras (rather than to
that of $C^*$-algebras) that cyclic homology applies.

The following result plays an important role in the theory of
endomotives and their applications to examples arising from the
study of the thermodynamical properties of certain quantum
statistical dynamical systems. We shall refer to the following proposition, in
section~\ref{BCsystem} of this paper for the description of the
properties of the ``BC-system''. The BC-system is a particularly relevant
quantum statistical dynamical system which has been the prototype
and the motivating example for the introduction of the notion of an
endomotive. We refer to \cite{7}, \S~4.1 for the definition of the
notion and the properties of a state on a (unital) involutive
algebra.

\begin{proposition}[\cite{7}~Proposition~3.1]\label{fab}\vspace{.1in}

1) The action \eqref{charaction} of $G_k$ on $X(\bar k)$ defines a
canonical action of $G_k$ by automorphisms of the $C^*$-algebra
$\bar\cA_\C=C(X)\rtimes S$, preserving globally $C(X)$ and such
that, for any pure state $\varphi$ of $C(X)$,
\begin{equation}\label{fabulous0}
\alpha \,\varphi(a)=\, \varphi (\alpha^{-1}(a)) \qqq a\in A \,,\quad
\alpha \in G_k.
\end{equation}

2) When the Artin motives $A_\alpha$ are abelian and normal, the
subalgebras $A\subset C(X)$ and $\cA_k\subset \bar\cA_\C$ are
globally invariant under the action of $G_k$ and the states
$\varphi$ of $\bar\cA_\C$ induced by pure states of $C(X)$ fulfill
\begin{equation}\label{fabulous}
\alpha \,\varphi(a)=\, \varphi (\theta(\alpha)(a)) \qqq a\in \cA_k
\,,\quad \theta(\alpha)=\alpha^{-1} \qqq \alpha \in
G^{ab}_k=\,G_k/[G_k,G_k]
\end{equation}
\end{proposition}

On the totally disconnected compact space $X$, the abelian semigroup
$S$ of homeomorphisms acts, producing closed and open subsets
$X^s\hookrightarrow X$, $x\mapsto x\cdot s$.  The normalized
counting measures $\mu_\alpha$ on $X_\alpha$ define a {\it
probability measure} on $X$ with the property that the
Radon--Nikodym derivatives
\begin{equation}\label{RNder}
\frac{ds^*\mu}{d\mu}
\end{equation}
are locally constant functions on $X$. One lets $\cG=X \rtimes S$ be
the corresponding \'etale locally compact groupoid. The crossed
product $C^*$-algebra $C(X(\bar k))\rtimes S$ coincides with the
$C^*$-algebra $C^*(\cG)$ of the groupoid $\mathcal G$.

The notion of right (or left) action of $\cG$ on a  totally
disconnected locally compact space $\cZ$ is defined as in the
algebraic case by \eqref{groupoidact11} and \eqref{groupoidact12}. A
right action of $\cG$ on $\cZ$ gives on the space $C_c(\cZ)$ of
continuous functions with compact support on $\cZ$ a structure of
right module over $C_c(\cG)$.
 When the fibers of the map $g: \cZ \to X$ are discrete
(countable) subsets of $\cZ$ one can define on $C_c(\cZ)$ an inner
product with values in $C_c(\cG)$ by
\begin{equation}\label{groupoidact3loc}
\langle \xi,\eta\rangle(x,s)=\,\sum_{z\in g^{-1}(x)}\,\bar
\xi(z)\,\eta(z\circ s)
\end{equation}

A right action of $\cG$ on $\cZ$ is {\em \'etale} if and only if the
fibers of the map $g$ are discrete and the identity is a compact
operator in the right $C^*$-module $\cE_\cZ$ over $C^*(\cG)$ given
by \eqref{groupoidact3loc}.

An {\it \'etale correspondence} is a
$(\cG(X_\alpha,S),\cG(X'_{\alpha'},S'))$-space $\cZ$ such that the
right action of $\cG(X'_{\alpha'},S')$ is \'etale.

The $\Q$-vector space
$$ \Corr((X,S,\mu),(X',S',\mu')) $$
of linear combinations of \'etale   correspondences $\cZ$ modulo the
equivalence relation $\cZ\cup \cZ'=\cZ+\cZ'$ for disjoint unions,
defines the space of (virtual) correspondences.

For $M=(X,S,\mu)$, $M'=(X',S',\mu')$, and $M''=(X'',S'',\mu'')$, the
composition of correspondences
$$\Corr(M,M')\times \Corr(M', M'') \to \Corr(M,M''),\qquad (Z,W)\mapsto Z\circ  W $$ is given
following the same rule as for the algebraic case \eqref{comprule},
that is by the fiber product over the groupoid $\cG'$. A
correspondence gives rise to a bimodule $\cM_\cZ$ over the algebras
$C(X)\rtimes S$ and $C(X')\rtimes S'$ and the composition of
correspondences translates into the tensor product of bimodules.

\begin{definition}\label{ExtArtCalg}
The category $C^*\cV^o_K$ of {\em analytic endomotives} is the
(pseudo)abelian category generated by objects of the form
$M=(X,S,\mu)$ with the properties listed above and morphisms given
as follows. For $M=(X,S,\mu)$ and $M'=(X',S',\mu')$ objects in the
category, one sets
\begin{equation}\label{CorrMCalg}
\Hom_{C^*\cV^0_K}(M,M')= \Corr(M,M')\otimes_\Q K.
\end{equation}
\end{definition}

The following result establishes a precise relation between the
categories of Artin motives and (noncommutative) endomotives.

\begin{theorem}[\cite{7}, Theorem 3.13]\label{NCArtinThm}
The categories of Artin motives and algebraic and analytic
endomotives are related as follows.
\begin{enumerate}
\item The map $\cG\mapsto \cG(\bar k)$ determines a
tensor functor
\begin{equation*}
\cF :\cE\cV^o(k)_{K} \to C^*\cV^o_K,\qquad \mathcal
F(X_\alpha,S)=(X(\bar k),S,\mu)
\end{equation*}
from algebraic to analytic endomotives.
\item The Galois group $G_k=\Gal(\bar k/k)$
acts by natural transformations of $\cF$.
\item The category $\cC\cV^o(k)_{K}$ of Artin motives embeds as a
full subcategory of $\cE\cV^o(k)_{K}$.
\item The composite functor
\begin{equation}\label{kFfunctor}
c\circ \cF: \cE\cV^o(k)_{K} \to \cK\cK\otimes K
\end{equation}
maps the full subcategory $\cC\cV^o(k)_{K}$ of Artin motives
faithfully to the category $\cK\cK_{G_k}\otimes K$ of
$G_k$-equivariant $KK$-theory with coefficients in $K$.
\end{enumerate}
\end{theorem}

 Given two Artin motives $X = \text{Spec}(A)$ and $X' =
\text{Spec}(B)$ and a component $Z = \text{Spec}(C)$ of the
cartesian product $X\times X'$, the two projections turn $C$ into a
$(A,B)$-bimodule $c(Z)$. If $U = \sum a_i\chi_{Z_i}\in
\text{Hom}_{\cC\cV^o(k)_{K}}(X,X')$, $c(U) = \sum a_i c(Z_i)$
defines a sum of bimodules in $\cK\cK\otimes K$. The composition of
correspondences in $\cC\cV^o(k)_{K}$ translates into the tensor
product of bimodules in $\cK\cK_{G_k}\otimes K$
\[
c(U)\otimes_B c(L) \simeq c(U\circ L).
\]
One composes the functor $c$ with the natural functor $A \to A_\C$
which associates to a $(A,B)$-bimodule $\mathcal E$ the
$(A_\C,B_\C)$-bimodule $\mathcal E_\C$. The resulting functor
\[
c\circ\cF\circ\iota: \cC\cV^o(k)_{K} \to \cK\cK_{G_k}\otimes K
\]
is faithful since a correspondence such as $U$ is uniquely
determined by the corresponding map of $K$-theory $K_0(A_\C)\otimes
K \to K_0(B_\C)\otimes K$.

\subsection{The endomotive of the BC system}\label{BCsystem}

The prototype example of the data which define an analytic
endomotive is the system introduced by Bost and Connes in \cite{6}.
The evolution of this $C^*$-dynamical system encodes in its group of
symmetries the arithmetic of the maximal abelian extension of
$k=\Q$.

\no This quantum statistical dynamical  system is described by the datum
 given by a noncommutative $C^*$-algebra of
{\it observables} $\bar\cA_\C = C^*(\Q/\Z)\rtimes_\alpha \N^\times$
and by the {\it time evolution} which is assigned in terms of a
one-parameter family of automorphisms $\sigma_t$ of the algebra. The
action of the (multiplicative) semigroup $S = \N^\times$ on the
commutative algebra $C^*(\Q/\Z)\simeq C(\hat\Z)$ is defined by
\[
\alpha_n(f)(\rho) = \begin{cases} f(n^{-1}\rho), &\text{if $\rho\in
n\hat\Z$}\\0, &\text{otherwise}\end{cases},\qquad\qquad
\rho\in\hat\Z = \varprojlim_n \Z/n\Z.
\]
For the definition of the associate endomotive, one considers the
projective system $\{X_n\}_{n\in\N}$ of zero-dimensional algebraic
varieties $X_n=\Spec(A_n)$, where $A_n=\Q[\Z/n\Z]$ is the group ring
of the abelian group $\Z/n\Z$. The inductive limit
$A=\displaystyle{\varinjlim_n} A_n = \Q[\Q/\Z]$ is the group ring of
$\Q/\Z$. The endomorphism $\rho_n: A \to A$ associated to an element
$n\in S$ is given on the canonical basis $e_r \in \Q[\Q/\Z]$, $r\in
\Q/\Z$, by
\begin{equation} \label{rhoBC}
\rho_n(e_r)=\,\frac{1}{n}\,\sum_{ns=r}\,e_s.
\end{equation}
The Artin motives $X_n$ are normal and abelian, so that
Proposition~\ref{fab} applies. The action of the Galois group $G_\Q
= \text{Gal}(\bar\Q/\Q)$ on $X_n=\Spec(A_n)$ is obtained by composing
a character $\chi: A_n \to \bar\Q$ with the action of an element
$g\in G_\Q$. Since $\chi$ is determined by the $n$-th root of unity
$\chi(e_{1/n})$, this implies that the action of $G_\Q$ factorizes
through the cyclotomic action and coincides with the {\it symmetry
group} of the BC-system. The subalgebra $\cA_\Q\subset \bar\cA_\C =
C(X)\rtimes S$ coincides with the rational subalgebra defined in
\cite{6}.

\no There is an interesting description of this system in terms of a
pointed algebraic variety $(Y,y_0)$ (\cf~section~\ref{genendo}) on
which the abelian semigroup $S$ acts by finite endomorphisms. One
considers the pointed affine group scheme $(\bG_m,1)$ (the
multiplicative group) and lets $S$ be the semigroup of non-zero
endomorphisms of $\bG_m$. These endomorphisms correspond to maps of
the form $u\mapsto u^n$, for some $n\in\N$. Then, the general
construction outlined in section~\ref{genendo} determines on
$(\bG_m(\Q),1)$ the BC system.

\no One considers the semigroup $S=\N^\times$ acting on $\bG_m(\Q)$ as
specified above.  It follows from the definition \eqref{Xs} that
$X_n =\Spec (A_n)$ where
\[
A_n = \Q[u_n^{\pm 1}]/(u_n^n-\,1).
\]
For $n|m$ the connecting morphism $\xi_{m,n}:X_m\to X_n$ is defined
by the algebra homomorphism $A_n \to A_m$ , $u_n^{\pm 1}\mapsto
u_m^{\pm a}$ with $a=\,m/n$. Thus, one obtains an isomorphism of
$\Q$-algebras
\begin{equation}\label{BClim}
\iota\;:\;A=\varinjlim_n \,A_n \stackrel{\sim}{\to}
\,\Q[\Q/\Z]\,,\quad \iota(u_n)=\,e_{\frac{1}{n}} .
\end{equation}
The partial isomorphisms $\rho_n: \Q[\Q/\Z] \to \Q[\Q/\Z]$ of the
group ring as described by the formula \eqref{rhoBC} correspond
under the isomorphism $\iota$, to those given by \eqref{endoS} on
$X=\displaystyle{\varprojlim_n}X_n$. One identifies $X$ with its
space of characters $\Q[\Q/\Z]\to \bar \Q$. Then, the projection
$\xi_m(x)$ is  given by the restriction of (the character associated
to) $x\in X$ to the subalgebra $A_m$. The projection of the
composite of the endomorphism $\rho_n$ of \eqref{rhoBC} with $x\in
X$ is given by
$$
x(\rho_n(e_r))=\,\frac{1}{n}\,\sum_{ns=r}\,x(e_s)\,.
$$
This projection is non-zero if and only if the restriction
$x|_{A_n}$ is the trivial character, that is if and only if
$\xi_n(x)=1$. Moreover, in that case one has
$$
x(\rho_n(e_r))=\,x(e_s)\qqq s \,,\quad ns = r\,,
$$
and in particular \begin{equation}\label{actrho}
x(\rho_n(e_{\frac{1}{k}}))=\,x(e_{\frac{1}{nk}}).
\end{equation}

For $k|m$ the inclusion of algebraic spaces $X_k\subset X_m$ is
given at the algebra level  by the surjective homomorphism
$$
j_{k,m}\;:\; A_m\to A_k\,,\quad j_{k,m}(u_m)=\, u_k .
$$
Thus, one can rewrite \eqref{actrho} as
\begin{equation}\label{actrho1}
x\circ\rho_n\circ j_{k,nk}=\,x|_{A_{nk}}.
\end{equation}
This means that
$$
 \xi_{nk}(x)=\, \xi_k(x\circ\rho_n).
$$
By using the formula \eqref{invrhos}, one obtains the desired
equality of the $\rho$'s of \eqref{rhoBC} and \eqref{endoS1}.

\no This construction continues to make sense for the affine algebraic
variety $\bG_m(k)$ for any field $k$, including the case of a field
of positive characteristic. In this case one obtains new systems,
different from the BC system.

\section{Applications: the geometry of the space of ad\`eles classes}

The functor
\begin{equation*}
\cF :\cE\cV^o(k)_{K} \to C^*\cV^o_K
\end{equation*}
which connects the categories of algebraic and analytic endomotives
establishes a significant bridge between the commutative world of
Artin motives and that of noncommutative geometry. When one moves
from commutative to noncommutative algebras, important new tools of
thermodynamical nature become available. One of the most relevant
techniques (for number-theoretical applications) is supplied by the
theory of Tomita and Takesaki for von Neumann algebras (\cite{36}) which
associates to a suitable state $\varphi$ (\ie a faithful weight) on
a von Neumann algebra $M$, a one-parameter group of automorphisms of
$M$ (\ie the modular automorphisms group)
\[
\sigma_t^\varphi: \R \to \text{Aut}(M),\qquad \sigma_t^\varphi(x) =
\Delta_\varphi^{it}x\Delta_\varphi^{-it}.
\]
$\Delta_\varphi$ is the modular operator which acts on the
completion $L^2(M,\varphi)$ of $\{x\in M:\varphi(x^*x)<\infty\}$,
for the scalar product $\langle x,y\rangle = \varphi(y^*x)$.

\no This general theory  applies in particular to the unital involutive
algebras $\mathcal A = C^\infty(X)\rtimes_{alg} S$ of
\eqref{densesub} and to the related $C^*$-algebras $\bar\cA_\C$ which
are naturally associated to an endomotive.

\no A remarkable result proved by Connes in the theory and classifications of factors
(\cite{8}) states that, modulo inner automorphisms of $M$, the
one-parameter family $\sigma_t^\varphi$ is {\it independent} of the
choice of the state $\varphi$. This way, one obtains a {\it
canonically defined} one parameter group of automorphisms classes
\[
\delta: \R \to \text{Out}(M)=\text{Aut}(M)/\text{Inn}(M).
\]
In turn, this result implies that the crossed product dual algebra
\[
\hat M = M \rtimes_{\sigma_t^\varphi}\R
\]
and the dual scaling action
\begin{equation}\label{dualaction}
\theta_\lambda: \R^*_+ \to \text{Aut}(\hat M)
\end{equation}
are {\it independent} of the choice of (the weight) $\varphi$.

\no When these results are applied to the analytic endomotive $\mathcal
F(X_\alpha,S)$ associated to an algebraic endomotive
$M=(X_\alpha,S)$, the above dual representation of $\R^*_+$ combines
with the representation of the absolute Galois group $G_k$. In the
particular case of the endomotive associated to the BC-system
(\cf~section~\ref{BCsystem}), the resulting representation of
$G_\Q\times\R^*_+$ on the cyclic homology $HC_0$ of a suitable
$\Q(\Lambda)$-module $D(\mathcal A,\varphi)$ associated to the
thermodynamical dynamics of the system $(\cA,\sigma_t^\varphi)$
determines the spectral realization of the zeroes of the Riemann
zeta-function and of the Artin $L$-functions for abelian characters
of $G_k$ (\cf\cite{7}, Theorem~4.16).

\no The action of the group $W=G_\Q\times\R^*_+$ on the cyclic homology
$HC_0(D(\mathcal A,\varphi))$ of the noncommutative motive
$D(\mathcal A,\varphi)$ is analogous to the action of the Weil group
on the \'etale cohomology of an algebraic variety. In particular,
the action of $\R^*_+$ is the `characteristic zero' analog of the
action of the (geometric) Frobenius on \'etale cohomology. This
construction determines a functor
\[
\omega: \cE\cV^o(k)_{K} \to \text{Rep}_\C(W)
\]
from the category of endomotives to the category of
(infinite-dimensional) representations of the group $W$.

\no The analogy with the Tannakian formalism of classical motive theory is
striking. It is also important to underline the fact that the whole
thermodynamical construction is non-trivial and relevant for
number-theoretic applications only because of the particular nature of the factor
$M$ (type $\text{III}_1$) associated to the original datum $(\cA,
\varphi)$ of the BC-system.

\no It is tempting to compare the original choice of the state $\varphi$
(weight) on the algebra $\mathcal A$ which singles out (via the
Gelfand-Naimark-Segal construction) the factor $M$ defined as the
weak closure of the action of $\bar\cA = C(X)\rtimes S$ in the
Hilbert space $\mathcal H_\varphi = L^2(M,\varphi)$, with the
assignment of a factor $(X,p,m)^r$ ($r\in\Z$) on a pure motive $M =
(X,p,m)$: \cf~\eqref{grading}. In classical motive theory, one knows
that the assignment of a $\Z$-grading is canonical only for
homological equivalence or under the assumption of the Standard
Conjecture of K\"unneth type. In fact, the definition of a weight
structure depends upon the definition of a complete system of
orthogonal central idempotents $\pi_X^i$.

\no Passing from the factor $M$ to the canonical dual representation
\eqref{dualaction} carries also the advantage to work in a setting
where projectors are classified by their real dimension ($\hat M$ is
of type $\text{II}_1$), namely in a noncommutative framework of
continuous geometry which generalizes and yet still retains some relevant
properties of the algebraic correspondences (\ie degree or
dimension).

\no The process of dualization is in fact subsequent to a
thermodynamical ``cooling procedure'' in order to work with a system
whose algebra approaches and becomes in the limit, a commutative
algebra (\ie $\text{I}_\infty$). Finally, one has also to implement
a further step in which one filters (\ie ``distils'') the relevant
noncommutative motive $D(\mathcal A,\varphi)$ within the derived
framework of cyclic modules (\cf~section~\ref{cyclicmodules}). This
procedure is somewhat reminissent of the construction of the
vanishing cohomology in algebraic geometry
(\cf\cite{23}).

\no When the algebra of the BC-system gets replaced by the
noncommutative algebra of coordinates $\mathcal A = \mathcal
S(\mathbb A_k)\rtimes k^*$ of the ad\`ele class space $X_k= \mathbb
A_k/k^*$ of a number field $k$, the cooling procedure is described
by a restriction morphism of (rapidly decaying) functions on $X_k$
to functions on the ``cooled down'' subspace $C_k$ of id\`ele
classes (\cf\cite{7}, Section~5). In this context, the
representation of $C_k$ on the cyclic homology $HC_0(\mathcal
H^1_{k,\C})$ of a suitable noncommutative motive $\mathcal
H^1_{k,\C}$  produces the spectral realization of the zeroes of
Hecke $L$-functions (\cf~{\it op.cit}, Theorem~5.6).

\no The whole construction describes also a natural way to associate to
a noncommutative space a canonical set of ``classical points'' which
represents the analogue in characteristic zero, of the geometric
points $C(\bar{\bf F}_q)$ of a smooth, projective curve $C_{/{\bf
F}_q}$.

\section{Bibliography}


\end{document}